\UseRawInputEncoding
\documentclass[12pt]{amsart}

\usepackage[all]{xy}
\usepackage{fullpage}
\usepackage{latexsym}
\usepackage{amsmath}
\usepackage{amsfonts}
\usepackage{amssymb}
\usepackage{amsthm}
\usepackage{eucal}
\usepackage{enumerate,yfonts}
\usepackage{mathrsfs}
\usepackage{graphicx}
\usepackage{graphics}
\usepackage{epstopdf}
\usepackage{amscd}
\usepackage{bbm}
\usepackage{hyperref}
\usepackage{url}
\usepackage{color}
\usepackage{bbm}
\usepackage{cancel}
\usepackage{enumerate}
\usepackage{amsmath,amsthm}
\usepackage{amssymb}
\usepackage{epsfig}
\usepackage{pstricks}
\usepackage{xy}
\usepackage{xypic}
\usepackage{bbm}
\usepackage{inputenc}

\newtheorem{thm}{Theorem}[section]
\newtheorem{corollary}[thm]{Corollary}
\newtheorem{lemma}[thm]{Lemma}

\newtheorem{proposition}[thm]{Proposition}
\newtheorem{prop}[thm]{Proposition}
\newtheorem{conjecture}[thm]{Conjecture}
\newtheorem{thm-dfn}[thm]{Theorem-Definition}

\theoremstyle{definition}
\newtheorem{definition}[thm]{Definition}

\numberwithin{equation}{section}

\theoremstyle{remark}
\newtheorem{remark}{Remark}[section]
\newtheorem{example}[remark]{Example}

\newcommand{\fg}{{\mathfrak g}}

\newcommand{\fp}{{\mathfrak p}}

\newcommand{\rW}{{\mathrm W}}

\newcommand{\bZ}{{\mathbb Z}}

\newcommand{\mE}{\mathcal{E}}
\newcommand{\mF}{\mathcal{F}}

\newcommand{\mM}{\mathcal{M}}

\newcommand{\mO}{\mathcal{O}}
\newcommand{\mL}{\mathcal{L}}

\newcommand{\on}{\operatorname}

\newcommand{\ra}{\rightarrow}
\newcommand{\la}{\leftarrow}

\newcommand{\is}{\simeq}

\newcommand{\Loc}{\on{Loc}}

\newcommand{\quash}[1]{}  %%Anything in \quash is ignored
\newcommand{\nc}{\newcommand}

\newcommand{\bfI}{\textbf{I}}
\newcommand{\bfP}{\textbf{P}}
\newcommand{\bfQ}{\textbf{Q}}

\newcommand{\frakA}{{\mathfrak A}}

\newcommand{\frakM}{{\mathfrak M}}

\newcommand{\frakZ}{{\mathfrak Z}}

\newcommand{\bbF}{{\mathbb F}}

\newcommand{\bbQ}{{\mathbb Q}}

\newcommand{\bbZ}{{\mathbb Z}}
\newcommand{\calA}{{\mathcal A}}

\newcommand{\calC}{{\mathcal C}}

\newcommand{\calF}{{\mathcal F}}

\newcommand{\calM}{{\mathcal M}}
\newcommand{\calN}{{\mathcal N}}

\newcommand{\calS}{{\mathcal S}}

\newcommand{\calZ}{{\mathcal Z}}

\nc{\al}{{\alpha}} \nc{\be}{{\beta}} \nc{\ga}{{\gamma}}
\nc{\ve}{{\varepsilon}} \nc{\Ga}{{\Gamma}} %\nc{\la}{{\lambda}}
\nc{\La}{{\Lambda}}
\nc{\ho}{\text{hol}}

\nc{\ad }{{\on{ad }}}
\nc{\Par }{{\on{Par }}}

\nc{\aff}{{\on{aff}}} \nc{\Aff}{{\mathbf{Aff}}}

\nc{\der}{{\on{der}}}

\nc{\diag}{{\on{diag}}}
\newcommand{\End}{{\on{End}}}
\nc{\Fl}{{\calF\ell}}

\nc{\Hg}{{\on{Higgs}}}
\newcommand{\Hom}{{\on{Hom}}}
\newcommand{\id}{{\on{id}}}
\nc{\Id}{{\on{Id}}}

\nc{\Ind}{{\on{Ind}}}

\nc{\Op}{{\on{Op}}}

\newcommand{\pr}{{\on{pr}}}
\newcommand{\Res}{{\on{Res}}}
\nc{\res}{{\on{res}}}

\nc{\tr}{{\on{tr}}}
\newcommand{\Tr}{{\on{Tr}}}

\newcommand{\GL}{{\on{GL}}}
\nc{\GSp}{{\on{GSp}}} \nc{\GU}{{\on{GU}}} \nc{\SL}{{\on{SL}}}
\nc{\SU}{{\on{SU}}} \nc{\SO}{{\on{SO}}}

\nc{\nh}{{\Loc_{J^p}(\tau')}}
\nc{\bnh}{{\Loc_{\breve J^p}(\tau')}}

\nc{\bU}{{\overline{U}}} \nc{\IC}{{\on{IC}}}

\newcommand{\Fr}{{\on{Fr}}}

\newcommand{\beqn}{\begin{equation*}}
\newcommand{\eeqn}{\end{equation*}}

\newcommand{\beq}{\begin{equation}}
\newcommand{\eeq}{\end{equation}}

\begin{document}
\title{Towards the depth zero stable Bernstein center conjecture}
 \author{Tsao-Hsien Chen}
        \address{School of Mathematics,
      University of Minnesota, Twin Cities
          }
        \email{chenth@umn.edu}

\maketitle   
\begin{abstract}
Let $G$ be a split connected reductive group over 
 a non-archimedan local field $F$.
The depth zero stable Bernstein conjecture asserts that 
there is an algebra isomorphism 
between the depth zero stable Bernstein center of $G(F)$
and the ring of functions 
on the moduli of tame Langlands parameters.  An approach to the depth zero stable Bernstein conjecture was proposed in the work of Bezrukavnikov-Kazhdan-Varshavsky \cite{BKV}.
In this paper we generalize  results and techniques  in \cite{BKV} 
and apply them to give a geometric construction of elements 
in the depth zero Bernstein center. We conjecture that 
 our construction produces all elements in the depth zero  Bernstein center.
 An an illustration of the method, we give a construction  
 of an algebra embedding from the (limit of) stable Bernstein centers for
finite reductive groups to the depth zero Bernstein center 
and  a family of elements in the depth zero Bernstein center coming from 
Deligne's epsilon factors.
The paper is the first step toward the 
depth zero stable Bernstein center.

\end{abstract} 
\setcounter{tocdepth}{2}
\tableofcontents

\section{Introduction}

\subsection{The stable Bernstein center conjecture}
Let $G$ be a split reductive group over a non-archimedan local field $F$.
Let $Z(G(F))$ be the Bernstein center of $G(F)$
and let $Z^{st}(G(G))\subset Z(G(F))$ be the subspace consisting of element $z\in Z(G(F))$ 
such that the associated invariant distribution $\nu_z$ on $G(F)$ is stable.
A version of the stable Bernstein center conjecture 
asserts that there exists an algebra isomorphism  
\beq\label{stable center}
\mO(\Loc_{\hat G,F})\is Z^{st}(G(F))
\eeq
from the ring of functions on the moduli stack of local Langlands parameters $\on{Loc}_{\hat G,F}$
to the stable Bernstein center $Z^{st}(G(F))$ (see, e.g., \cite{BKV,FS,H,SS, Z2}).
Let $Z^0(G(F))\subset Z(G(F))$ be the subalgebra of depth zero Bernstein center. 
The moduli stack $\Loc^t_{\hat G,F}$ of tame Langlands parameters 
is a component of $\Loc_{\hat G,F}$ and 
it is expected that 
the isomorphism~\eqref{stable center} restrict to
an isomorphism
\beq\label{depth zero stable center}
\mO(\Loc^t_{\hat G,F})\simeq Z^{st,0}(G(F))
\eeq
where
$Z^{st,0}(G(F))=Z^{0}(G(F))\cap Z^{st}(G(F))$.
We will refer to the isomorphism~\eqref{depth zero stable center} the depth zero Bernstein center conjecture.\footnote{Note that  the stable Bernstein center conjecture  implies that the subspaces 
$Z^{st}(G(F))$ and $Z^{st,0}(G(F))$ are unital subalgebras of $Z(G(F))$.
The later assertion is the version of stable Bernstein center conjecture in \cite{BKV}.
}

An approach to the depth zero Bernstein center conjecture using $\ell$-adic sheaves 
was proposed in
the work of Bezrukavnikov-Kazhdan-Varshavsky  \cite{BKV}. 
In this paper we generalize various results and techniques  in \cite{BKV} 
and uses them to give a geometric construction of 
 elements in $Z^0(G(F))$. We conjecture that 
 our construction produces all elements in $Z^0(G(F))$.
 An an illustration of the method, we give a geometric construction 
 of an embedding from the (limit) stable Bernstein centers for
finite reductive groups to $Z^0(G)$
and  a family of elements in depth zero Bernstein center coming from 
Deligne's epsilon factors.
In the sequel \cite{C2,C3}, we will use the results in the paper 
to construct an algebra map
$\mO(\Loc^t_{\hat G,F})\to Z^0(G(F))$ and verify some properties 
predicted by the depth zero local Langlands correspondence.

The paper and its sequels are the first step toward to depth zero stable Bernstein center conjecture. 
We now describe the paper in more details.

 \subsection{Main results}
 In this paper we assume 
$G$ is a split, semisimple, and simply connected group over $\bbF_q$. 
Let $LG$ and
$L^+G$ be the loop group and formal arc group of $G$.
Let $\bfI\subset L^+G$ be an Iwahori subgroup and 
let $\bfI^+$ be its pro-unipotent radical.
Let $F=\bbF_q((t))$ and let 
$G(F)=LG(k)$ be the corresponding 
reductive group over $F$.
In their work on  stable center conjecture \cite{BKV}, Bezrukavnikov-Kazhdan-Varshavsky 
introduced a version of categorical center of the (universal) affine Hecke category, denoted by 
$\calZ_{\bfI^+}(LG)$, and they 
outlined a construction of an algebra homomorphism
\beq\label{[Z]}
K_0(\calZ_{\bfI^+}^\Fr(LG))\to Z^0(G(F))
\eeq
where $K_0(\calZ_{\bfI^+}^\Fr(LG))$ is the Grothendieck group (tensored over $\overline\bbQ_\ell$) of the category of 
Frobenius equivariant objects  of 
$\calZ_{\bfI^+}(LG)$. %\footnote{The definition of categorical center $\calZ_{\bfI^+}(LG)$in \emph{loc. cit.} is quite involved, for example, it is crucial to use the language of $\infty$-categories.}
They conjectured that the map~\eqref{[Z]} is surjective  %\cite[Conjecture 3.4.8]{BKV} 
and hence provides a geometric construction 
of all elements in $Z^0(G(F))$.
The construction of~\eqref{[Z]} relies on a conjectural stabilization theorem\footnote{
The conjectural stabilization theorem was stated in \cite["Theorem 3.3.16``]{BKV}, where ''Theorem`` means work in progress.}  for objects in the categorical center $\calZ_{\bfI^+}(LG)$ %, see \cite[''Theorem 3.3.16``]{BKV}\footnote{In \cite{BKV}, ''Theorem`` means work in progress}, 
and 
one of the main result in their paper is a proof of a 
Grothendieck group version of the stabilization theorem
for the monoidal unit element in
$\calZ_{\bfI^+}(LG)$ (\cite[Theorem 4.1.5]{BKV}). As an application, they
gave a geometric construction of the 
 Bernstein projector to the depth zero spectrum and verify its stability
 (\cite[Theorem 4.4.1]{BKV}).

 Inspired by the work of Bezrukavnikov-Kazhdan-Varshavsky,
 in this paper 
we introduce 
and study the following algebra:
\[\frak A(LG)=\on{lim}_{\bfP\in\Par} K_0(M(\frac{LG/\bfP^+}{\bfP}))\]
where 
$K_0(M(\frac{LG/\bfP^+}{\bfP}))$ is the Grothendieck group (tensored over $\overline\bbQ_\ell$) of
the parahoric Hecke category $M(\frac{LG/\bfP^+}{\bfP})$ 
associated to
each standard parahoric subgroup $\bfP\in\Par$ and the limit is taken with respect to
the so called Harish-Chandra transforms 
$\on{HC}_{\bfP,\bfQ}$,
 where $\bfP\subset\bfQ\in\Par$
(see Section \ref{affine HC}). 
%Objects in $A(LG)$ are called  \emph{admissible complexes}.
 We observe that the proof of the stabilization theorem for the unit element 
in \cite[Theorem 4.1.1]{BKV}
 can be applied to
 objects in $\mathfrak A(LG)$
 and we prove the following generalization:
  Let $(M(LG),*)$ be the monoidal category of constuctible $\overline\bbQ_\ell$-sheaves 
on $LG$ and let $K_0(M(LG))$ be its Grothendieck group. 

\begin{thm}[Theorem \ref{stabilization}]
\label{main result 1}
Let $\mM=\{\langle\mM_\mathbf P\rangle\}_{\mathbf P\in\Par}\in\frak A(LG)$.
For any $\mF\in M(LG)$, 
the system $\{\langle A^\mathrm Y_\mM\rangle*\langle\mF\rangle\}_{\mathrm Y\in\Upsilon}$
 of objects in $K_0(M(LG))$
stabilizes.
Here $\Upsilon$ is the partially ordered set of closed $\mathbf I$-invariant sub-schemes
of the affine flag $LG/\mathbf I$ and $\langle A^\mathrm Y_\mM\rangle\in K_0(M(LG))$ is the object in~\eqref{A^Y_M}.
\end{thm}

%Theorem \ref{main result 1} is restated as 
%Theorem \ref{stabilization} in the text; we direct reader to 
%Section \ref{Bernstein and admissible} for a more detailed statement and the proof.
Now using the above theorem one obtains the following construction 
of elements in the Grothendieck group version of the Bernstein centers.
Following \cite[Section 3.4.2]{BKV}, let us denote by 
$\frak Z(LG)=\End_{K_0(M(LG))^2}(K_0(M(LG)))$ the $\overline\bbQ_\ell$-algebra
of endomorphism of $K_0(M(LG))$, viewed as a $K_0(M(LG))^2$-module.
For any $\mM=\{\langle\mM_\mathbf P\rangle\}_{\mathbf P\in\Par}\in\frak A(LG)$,
we define $\langle A_\mM\rangle\in\on{End}(K_0(M(LG)))$
by the formula
\[\langle A_\mM\rangle(\langle\mF\rangle):=\lim_{Y\in\Upsilon^{op}}\langle A_\mM^Y\rangle*\langle\mF\rangle.\]
Note that it is well-defined thanks to the stabilization theorem above.
The following theorem generalizes \cite[Theorem 4.1.9]{BKV}:
\begin{thm}[Theorem \ref{main}]\label{main result 2}
The element $\langle A_\mM\rangle\in\on{End}(K_0(M(LG)))$ belongs to $\frak Z(LG)$. Moreover,
 the assignment $\mM\to\langle A_\mM\rangle$
defines an algebra map
$\langle A\rangle:\frak A(LG)\to\frak Z(LG)$.
\end{thm}
%This is restated as Theorem \ref{main} in the text.
%We give a full explanation of the result in Section \ref{A to Z}.

Now by applying the sheaf-function correspondence, we obtain the following 
construction of depth zero Bernstein center: 
Denote by
$\frakA^\Fr(LG)=\on{lim}_{\bfP\in\Par} K_0(M^\Fr(\frac{LG/\bfP^+}{\bfP}))$
and $\frakZ^\Fr(LG)=\End_{(K_0(M^\Fr(LG)))^2}(K_0(M^\Fr(LG)))$
where $M^\Fr(\frac{LG/\bfP^+}{\bfP})$ and $M^\Fr(LG)$
are the category of 
Frobenius equivariant objects (i.e. Weil objects)
in $M(\frac{LG/\bfP^+}{\bfP})$ and $M(LG)$ respectively.
Then the map $\langle A\rangle$ in Theorem \ref{main result 2} admits a lift
\[\langle A^\Fr\rangle:\frakA^\Fr(LG)\to\frakZ^\Fr(LG).\]
For any $\bfP\in\Par$, let  $\mathrm P=\bfP(\bbF_q)$
and $\mathrm P^+=\bfP(\bbF_q)$ be the corresponding 
parahoric subgroup  and  pro-unipotent radical respectively. We denote by
$(M(\frac{G(F)/\mathrm P^+}{\mathrm P}),*)$ be the parahoric Heck algebra 
of $G(F)$ consisting of $\mathrm P^+$ bi-invaraint
and $\mathrm P$-conjugation invariant   smooth measures on $G(F)$ with compact support.
Consider the following algebra 
\[A(G(F))=\lim_{\bfP\in\Par} M(\frac{G(F)/\mathrm P^+}{\mathrm P})\]
where the limit is taken with respect to the map
$M(\frac{G(F)/\mathrm Q^+}{\mathrm Q})\to M(\frac{G(F)/\mathrm P^+}{\mathrm P})$ sending 
$h$ to $h*\delta_{\mathrm P^+}$, where $\bfP\subset\bfQ\in\Par$ and
$\delta_{\mathrm P^+}$ is the Haar measure of $\mathrm P^+$ with total measure one.
The following theorem generalizes \cite[Theorem 4.4.1]{BKV}:

\begin{thm}[Theorem \ref{elementary construction}]\label{main result 3}
The is an algebra map
$[A]:A(G(F))\to Z^0(G(F))$ fitting into the following commutative diagram
\beq\label{intro diagram 1}
\xymatrix{
\frakA^\Fr(LG)\ar[d]\ar[r]^{\langle A^\Fr\rangle}&\frakZ^\Fr(LG)\ar[d]\\
A(G(F))\ar[r]^{[A]}&Z^0(G(F))}
 \eeq
where the vertical arrows are given by 
the sheaf-function correspondence.
\end{thm}

In \emph{loc. cit.} we also give 
an explicit formula for $[A]$.
In view of \cite[Conjecture 2.4.8]{BKV},
we propose the following:
\begin{conjecture}
The composed map $\frakA^\Fr(LG)\stackrel{\langle A^\Fr\rangle}\to\frakZ^\Fr(LG)\to Z^0(G(F))$ in~\eqref{intro diagram 1}
is surjective.
\end{conjecture}

\begin{remark}
We do not known whether  the vertical maps in~\eqref{intro diagram 1} are surjective.
Thus the construction of the map $[A]$ is not a direct consequence of 
Theorem \ref{main result 2} and sheaves-functions correspondence.
Instead, we just simply mimick the proof of Theorem \ref{main result 2}
but omit all the geometry.
\end{remark}
\quash{
\begin{remark}
The algebras $\frakA(LG)$, $\frakA(LG)^\Fr$,
and $A(G(F))$ are priori non-commutative and hence 
the maps $\langle A\rangle$, $\langle A\rangle^\Fr$ and $[A]$ 
in Theorem \ref{main result 2} and Theorem \ref{main result 3}
factor through the abelianization of 
$\frakA(LG)$, etc.
\end{remark}
}
\begin{remark}
The relationship between Theorem \ref{main result 3}
and the construction of~\eqref{[Z]} is as follows.
In \cite[Section 3.4.7]{BKV}, it was shown that 
the map~\eqref{[Z]} comes from a map
\[\langle Z^\Fr\rangle:K_0(\calZ^\Fr_{\bfI^+}(LG))\to\frak Z^\Fr(LG).\]
For any $\bfP\in\Par$,
we can form  the monoidal 
$\infty$-category $\calM(\frac{LG/\bfP^+}{\bfP})$ whose homotopy category is
$M(\frac{LG/\bfP^+}{\bfP})$
and also the 
monoidal 
$\infty$-category 
\beq\label{infty category version}
\calA(LG)=\lim_{\bfP\in\Par}\calM(\frac{LG/\bfP^+}{\bfP})
\eeq
where the limit is taken with respect to the 
Harish-Chandra transforms.\footnote{Since derived categories do not have limits in general,
in order define $\calA(LG)$ one has to use $\infty$-categories.} 
The monoidal $\infty$-category $\calA(LG)$ can be viewed as the categorification of 
the algebra
$\frakA(LG)$.
We have a natural map 
$K_0(\calA^\Fr(LG))\to\frakA^\Fr(LG)$ and, 
using \cite[``Theorem 3.3.9'']{BKV}, one can show that there is a natural monoidal
functor
\beq\label{functor S}
\mathcal Z_{\bfI^+}(LG)\to \calA(LG)
\eeq
such that 
$\langle Z^\Fr\rangle:K_0(\calZ^\Fr_{\bfI^+}(LG))\stackrel{~\eqref{functor S}}\to 
 K_0(\calA^\Fr(LG))\to
\frakA^\Fr(LG)\stackrel{\langle A^\Fr\rangle}\to \frak Z^\Fr(LG)$.
%We expect that the functor~\eqref{functor S} is an equivalence.
\end{remark}

Our motivation to introduce the algebras $\frak A(LG)$ 
and its function theoretical counterpart $A(G(F))$
comes from its connection with 
our previous work on Braverman-Kazhdan conjecture \cite{C1,C2}.
Namely, the algebra $\frak A(LG)$ and $A(G(F))$
contain the following subalgebras
\[\frak A(LG)_1=\{\{\langle\mM_\bfP\rangle\}_{\bfP\in\Par}\in\frak A(LG)|\on{supp}(\mM_{\bfP})\subset \frac{\bfP/\bfP^+}{\bfP}\text{\ for all\ }\mathbf P\in\Par\}\]
\[A(G(F))_1=\{\{ h_\bfP\}_{\bfP\in\Par}\in A(G(F))|\on{supp}(h_{\bfP})\subset\mathrm P\text{\ for all\ }\mathbf P\in\Par\}.\]
It turns out that 
the subalgebra $\frak A(LG)_1$ admits a description in terms of 
a class of sheaves on the reductive quotients of standard parahoric subgroups.
More precisely,  
for any reductive group $H$,
let
$D(\frac{H}{H})$
be the $H$-equivariant derived category of $H$ with respect to the 
adjoint action. Introduce the following subcategory of $A(H)\subset D(\frac{H}{H})$: 
\[A(H)=\{\mF\in D(\frac{H}{H})|\on{Supp}(\on{HC}_{P,H}(\mF))\subset\frac{P/U_{P}}{P}\text{\ for all\ } P\in\on{par}\}\]
%\beq\label{A(H)}
%A(H)=\{\mF\in M(\frac{H}{H})| \on{Supp}(\on{HC}_{P,H}(\mF))\subset\frac{P/U_P}{P} \ \text{for all standard parabolic subgroup}\ P \text{\ of\ }H\}
%\eeq
where $\on{HC}_{P,H}$ is the Harish-Chandra transform associated to  the
standard parabolic subgroup
$P\subset H$ (see Section \ref{HC fd}).
Objects in $A(H)$ 
 were appeared in the work of 
Braverman-Kazhdan \cite{BK1,BK2} on non-abelian Fourier transform 
for finite reductive groups and were studied in \cite{CN,C1,C2}.
It follows from the definition of $\frak A(LG)_1$ that 
there is an isomorphism of algebras:
\[\on{lim}_{\bfP\in\Par}K_0(A(L_\bfP))\is \frak A(LG)_1\]
here $L_\bfP=\bfP/\bfP^+$ is the reductive quotient of $\bfP$
and the limit is taken with respect to 
parabolic restriction functors (see Lemma \ref{1=par}). 
Using our previous work  \cite{C1,C2} on the Braverman-Kazhdan conjecture, 
we are able to produce 
many objects in $\on{lim}_{\bfP\in\Par}K_0(A(L_\bfP))$, and hence object
in $\frakA(LG)_1$,
in terms of a certain class of Weyl group equivariant complexes on the maximal 
torus  called the \emph{strongly central complexes}
(see Definition \ref{strongly central complexes}). Moreover, we 
show that under the function-sheaves correspondence  
the category of strongly central complexes
provides a 
categorification of the stable Bernstein center 
of the 
finite reductive group (see Proposition \ref{geometrization 1}). Now combining with 
Theorem \ref{main result 3}, we obtain the following geometric construction of a map
from the (limit of) stable Bernstein centers of finite reductive groups
to depth zero Bernstein center:
For any $\bfP\in\Par$, we denote by
$Z^{st}(L_\mathbf P(\bbF_q))$ the
stable Bernstein center for the finite reductive group $L_\mathbf P(\bbF_q)$
(see Section \ref{center for finite gps}).
Consider the following algebra
\[\lim_{\mathbf P\in\Par}Z^{st}(L_\mathbf P(\bbF_q))\]
where the limit is taken with respect to
the natural transfer map $\hat\rho_{\bfP,\bfQ}:Z^{st}(L_\mathbf Q(\bbF_q))\to Z^{st}(L_\mathbf P(\bbF_q))$,
$\mathbf P\subset\mathbf Q\in\Par$ (see~\eqref{transfer Q to P}).

\begin{thm}[Theorem \ref{decategorification}]\label{main result 4}
There is a natural injective algebra map
\beq\label{lim to Z}
\Psi:\lim_{\mathbf P\in\Par}Z^{st}(L_\mathbf P(\bbF_q))\to Z^0(G(F))
\eeq
characterized by the following formula:
for any element $z=\{z_\mathbf P\}\in\lim_{\mathbf P\in\Par}Z^{st}(L_\mathbf P(\bbF_q))$, 
 and a vector 
$v\in V$ in an irreducible representation $(\pi,V)$ of $G(F)$ of depth zero,
we have 
\beq\label{formula}
(\Psi(z)|_\pi)(v)=(z_\mathbf P|_{\pi^{\mathrm P^+}})(v)
\eeq
where $\mathbf P\in\Par$ is any parahoric subgroup such that 
$v\in V^{\mathrm P^+}$
and $\pi^{\mathrm P^+}$ denotes the natural representation of 
$L_\mathbf P(\bbF_q)$ on $V^{\mathrm P^+}$.
\quash{
We have the following commutative diagram 
\beq\label{key diagram}
\xymatrix{\lim_{\mathbf P\in\Par} 
K_0(A(T_\mathbf P/\rW_\mathbf P)^\Fr)\ar[r]^{}\ar[d]^{\Upsilon}& \frakA(LG)^\Fr\ar[d]\ar[r]&\frakZ(LG)^\Fr\ar[d]\\
\lim_{\mathbf P\in\Par}Z^{st}(L_\mathbf P(\bbF_q))\ar[r]^{}&A(G(F))\ar[r]&Z^0(G(F))}
 \eeq}
\end{thm}

%Theorem \ref{main result 4} is stated as Theorem \ref{decategorification}
%in the text; we direct the reader to Section \ref{SCC} for more detailed statement.
The natural projection map 
$\on{lim}_{\mathbf P\in\Par}Z^{st}(L_\mathbf P(\bbF_q))\to Z^{st}(G(\bbF_q))$
admits a section 
$Z^{st}(G(\bbF_q))\to \on{lim}_{\mathbf P\in\Par}Z^{st}(L_\mathbf P(\bbF_q))$ and,
by composing with the map~\eqref{lim to Z}, we obtain an embedding
\beq\label{Z^st to Z^0}
\zeta:Z^{st}(G(\bbF_q))\to Z^0(G(F))
\eeq
from the stable Bernstein center of $G(\bbF_q)$ to depth zero Bernstein center of $G(F)$.
The existence of the map~\eqref{Z^st to Z^0}
follows from the work of 
Moy-Prasad \cite{MP1,MP2} and Theorem \ref{main result 4} provides an alternative geometric construction.

As an application of Theorem \ref{main result 4} we give 
a geometric construction of  a family of elements in depth zero Bernstein center coming from 
Deligne's epsilon factors. Namely, 
let $\rho:\hat G\to\GL_n$ be a representation of the dual group and 
let $\mF_{G,\rho,\psi}\in D^\Fr(\frac{G}{G})$ be the corresponding Braverman-Kazhdan's Bessel sheaf
on $G$ introduced in \cite[Section 6]{BK2}.
It was shown in  \cite{C2,LL}  that the 
associated class function $\Tr(\Fr,\mF_{G,\rho,\psi}):G(\bbF_q)\to\overline\bbQ_\ell$
is stable (as conjectured by Braverman-Kazhdan) and hence gives rise to an element 
in the stable Bernstein center $z_{G,\rho}\in Z^{st}(G(\bbF_q))$ (see Section \ref{center for finite gps}). We denote by $z_\rho=\zeta(z_{G,\rho})\in Z^0(G(F))$
its image under the embedding~\eqref{Z^st to Z^0}.
On the other hand, 
using Deligne's epsilon factors $\epsilon_0(r,\psi_F,dx)$,  one can attach 
to each representation $\rho$ of the dual group $\hat G$ an element $z_{0,\rho}\in Z^0(G(F))$ (see Section \ref{Deligne}).
Let $\check\rho:\hat G\to\GL_n$ be the contragredient representation 
of $\rho$.

\begin{thm}[Theorem \ref{deligne}]\label{epsilon intro}
We have $z_\rho=(-1)^nz_{0,\check\rho}$.
\end{thm}

\begin{remark}
In the case when $G=\GL_n$ and $\rho=\id:\GL_n\to\GL_n$ is the standard representation, Theorem \ref{epsilon intro} is essentially
a reformulation of a result of Macdonald \cite{M}.
\end{remark}
Along the way, we complete the proof of 
a conjecture of Braverman-Kazhdan for arbitrary parabolic subgroups \cite[Conjecture 6.5]{BK2}:
\begin{thm}[Corollary \ref{BK conj}]
Let $P\subset G$ be a  parabolic subgroup with unipotent radical $U_P$
and let $f:G\to G/U_P$ be the quotient map.
The derived push-forward $f_!(\mF_{G,\rho,\psi})$ is supported on
$P/U_P\subset G/U_P$.

\end{thm}

\begin{remark}
In the case when $P=B$ is a Borel subgroup, this is the main result 
in \cite{C2}. In the case when  $G=\GL_n$ and arbitrary parabolic subgroup $P$,
 this is proved in \cite[Proposition A1]{CN}.
 The general case follows from the case when $P=B$ plus a simple 
 fact on Harish-Chandra transforms (see Corollary \ref{B implies P}).
\end{remark}

\subsection{Further directions} 
%In this final section of the introduction, 
We discuss results to
appear in sequel papers \cite{C3,C4} that build on those of the current paper
and also some possible extensions and generalizations.

%%%%%%%%
\quash{

\subsection{Applications}
We discuss some applications.
First we observe that 
the natural projection map 
$\on{lim}_{\mathbf P\in\Par}Z^{st}(L_\mathbf P(\bbF_q))\to Z^{st}(G(\bbF_q))$
admits a section 
$Z^{st}(G(\bbF_q))\to \on{lim}_{\mathbf P\in\Par}Z^{st}(L_\mathbf P(\bbF_q))$ and,
by composing with the map~\eqref{lim to Z}, we obtain an algebra map
\beq
Z^{st}(G(\bbF_q))\to Z^0(G(F))
\eeq
from the stable Bernstein center of $G(\bbF_q)$ to depth zero Bernstein center of $G(F)$.

The above  applications and their connection with 
moduli of Langlands parameters and
the so called Braverman-Kazhdan program 
are the original motivation for writing the papers.
}
%%%%%%

%\subsubsection{Gaitsgory's central sheaves}

\subsubsection{Semi-simple Langlands parameters}\label{moduli}

In the sequel
\cite{C3}, we will use the results in this paper to 
construct 
an algebra homomorphism 
\beq\label{phi}
\phi:\mO(\Loc^t_{\hat G,F})\to Z^{0}(G(F))
\eeq
fitting into the following commutative diagram
\[\xymatrix{\overline\bbQ_\ell[(\hat T//\rW)^{[q]}]\ar[d]\ar[r]^{\simeq}&Z^{st}(G(\bbF_q))\ar[d]^{\eqref{Z^st to Z^0}}\\
\mO(\Loc^t_{\hat G,F})\ar[r]^\phi& Z^{0}(G(F))}\]
where $(\hat T//\rW)^{[q]}$ is the set of semi-simple conjugacy classes in the dual group $\hat G$ (over $\overline\bbQ_\ell$) stable under $x\to x^q$,
the upper horizontal arrow is the natural isomorphism in~\eqref{Z=T//W}, and the 
left vertical arrow comes from the natural map 
\beq\label{projection}
\Loc^t_{\hat G,F}\to \Phi^t_{\hat G,F,I}\is(\hat T//\rW)^{[q]}
\eeq
sending 
a tame Langlands parameter $\rho$ to the inertia equivalence class 
of the semisimplification $r=\rho^{ss}$ (see Section \ref{Deligne}).
This is a first step toward the construction of the isomorphism~\eqref{depth zero stable center}.
As an application, we can associate to each  irreducible representation 
$\pi$ of $G(F)$ of depth zero a semi-simple Langlands parameter 
$\rho(\pi)\in \Loc^t_{\hat G,F}$ characterized by the properties that (1) the map 
$\mO(\Loc^t_{\hat G,F})\to Z^{0}(G(F))\to\End(\pi)=\overline\bbQ_\ell$
is given by evaluation at $\rho(\pi)$ and (2) the image 
of $\rho(\pi)$ along~\eqref{projection}  
is equal to the Deligne-Lusztig parameter $\theta(\pi)$ of $\pi$
(see Section \ref{DL parameters}).

\begin{remark}
Now
there are many  constructions of 
semi-simple Langlands parameters for depth zero representations:
the constructions by Fargues-Scholze \cite{FS}
and Lafforgue-Genestier \cite{LG} which work for any  irreducible representations of $G(F)$,
the construction by Lusztig \cite{Lu2} and 
 Hemo-Zhu \cite{HZ} for unipotent representations, and the one presented above.
The techniques used in those constructions are very different and it will be interesting to 
understand their relationship. 

\quash{
In the work \cite{FS,LG}, Fargues-Scholze and Lafforgue-Genestier
 attached to each  irreducible representations of $G(F)$
a semi-simple Langlands parameters. On the other hand,
the work of Lusztig \cite{Lu2} and the forthcoming work of Hemo-Zhu \cite{HZ}
on unipotent representations associated to each unipotent representation 
a semi-simple Langlands parameter.
We expect our semi-simple parameters $\rho(\theta)$
for depth zero representations  agree with  those constructed in \emph{loc. cit.}.}
\end{remark}

\subsubsection{Stability}
We conjecture that the image the map~\eqref{phi} lies in 
the subspace $Z^{st,0}(G(F))$ of stable Bernstein center.
In the sequel \cite{C4} we will show that the image of~\eqref{Z^st to Z^0} lies in $Z^{st,0}(G(F))$.
This generalizes (some of) the results in \cite{BKV,BV}.

\subsubsection{Mixed characteristic}
We expect that the results and techniques 
in this paper and the sequel \cite{C3}
can be extended
to the mixed characteristic case using the Witt vector version of
the Affine Grassmannian introduced by Zhu \cite{Z2}.
For example, the construction of the map 
$[A]:A(G(F))\to Z^0(G(F))$ is in fact valid for all  non-Archimedian local fields
(see Remark \ref{mixed char}).
However, the results in \cite{C4} on the  
stability of the image of~\eqref{Z^st to Z^0} 
 rely on the work of 
Yun \cite{Y} which is available (at the moment) only in the  equal characteristic case.\footnote{The proof in \cite{Y} uses the theory of Hitchin fibration which only exists  in 
the equal characteristic case.}

\subsection{Organization}
We briefly summarize here the main goals of each section.
In Section 2, we collect standard notation in loop groups.
In Section 3, we give a review of theory of $\ell$-adic 
sheaves on admissible ind-schemes and ind-stacks developed in \cite{BKV}
In Section 4, we prove the stabilization theorem for objects in the algebra $\frakA(LG)$.
In Section 5,
We give a geometric construction of elements in depth zero Bernstein centers.
In Section 6, we
introduce and
study  strongly central complexes 
and use them to produce 
objects in the algebra $\frakA(LG)$.
We show that the
category of strongly central complexes
provides a categorification of stable Bernstein centers for finite 
reductive groups and, 
combining with the results in Section 3, we construct a map from 
 stable Bernstein centers for finite 
reductive groups to depth zero Bernstein centers.
In Section 7,  
we give a geometric construction of  a family of elements in depth zero Bernstein center coming from 
Deligne's epsilon factors.

{\bf Acknowledgement.}
I would also like to 
thank the Institute of Mathematics Academia Sinica in Taipei for support, hospitality, and a nice research environment. 
A part of the work was done while the author visited the institute.
I am grateful to the participants of the
''Character sheaves on loop groups seminar`` 
in the Fall of 
2022 in Taipei.
Special thanks are due to Harrison Chen and Cheng-Chiang Tsai for 
organizing the seminar together and for 
several insightful discussions on classical and geometric Langlands correspondence. 
I also would like to thank David Nadler for useful discussion. 
The research is supported by NSF grant  DMS-2001257
and DMS-2143722.

\section{Notations}
\subsection{Loop groups}\label{loop groups}
Let $k$ be an algebraically closed field 
and
let $G$ be a semi-simple simply connected 
group over $k$. 
We assume the  characteristic $\on{char}(k)$ of $k$ is large and 
we fix a prime number $\ell$ different from $\on{char}(k)$.
We fixed a Borel subgroup $B$ and a maximal torus 
$T\subset B$. Let $\rW=N(T)/T$ be the Weyl group.
We write $\on{par}$ for the  partial ordered set of parabolic subgroups $P$ containing $B$.
For any $P\in\on{par}$ we write $U_P\subset P$ for the unipotent radial 
and denote by $L_P=P/U_P$ the reductive quotient.

Let $LG$ and $L^+G$ be the loop group and arc group of $G$.
Let $\textbf{I}\subset L^+G$ be the Iwahoric subgroup given by
 the preimage of the quotient map 
 $L^+G\to G$.
We denote by $\bfI^+$ the pro-unipotent radical of $\bfI$.
Let $\Delta$ be the set of simple roots associated to the 
pair $(T,B)$ and let $\widetilde\Delta$ be the set of affine simple roots.

We write $\on{Par}$ for the 
partial ordered 
set of parahoric subgroups $\bfP$ containing $\bfI$.
We write $\bfP^+$ for the pro-unipotent radical of 
$\bfP$ and denote by $L_{\bfP}=\bfP/\bfP^+$ the reductive quotient.
The quotient $B_\bfP=\bfI/\bfP^+\subset \bfP/\bfP^+=L_\bfP$  is a Borel subgroup of $L_\bfP$. The image of the composed map 
$L^+T\to\bfI\to B_\bfP $ 
is a maximal torus $T_\bfP\subset B_\bfP$ of $L_\bfP$.
We write $\rW_\bfP=N(T_\bfP)/T_\bfP$ the Weyl group of $L_\bfP$.
We set $\on{rk}(\bfP)=\on{rk}(L_\bfP^{\on{der}})$

Note that the 
natural surjection $L^+T\to T_\bfP$ 
factors through the quotient 
map $L^+T\to T$ and 
induces an isomorphism 
$T\is T_\bfP$. We write 
$\phi_\bfP:T_\bfP\is T$ for the inverse isomorphism.

Let $\Lambda$ be the weight lattice of $T$ and let
$\widetilde\rW=\rW\ltimes\Lambda$ be the affine Weyl group.
Note that 
for any  $\bfP\in\on{Par}$
the corresponding Weyl group $\rW_\bfP$ is naturally a subgroup of 
$\widetilde\rW$
and the identification $\phi_\bfP:T_\bfP\is T$ above is $\rW_\bfP$-equivariant 
where $\rW_\bfP$ acts on $T$ via the composed map
$\rW_\bfP\to\widetilde\rW\to\rW$ (where the last arrow is the projection map). 

More generally, for any $\bfP\subset\bfQ\in\Par$, 
the quotient $B_{\bfP,\bfQ}=\bfP/\bfQ^+\subset L_\bfQ=\bfQ/\bfQ^+$
is a parabolic subgroup 
with unipotent radical $U_{\bfP,\bfQ}=\bfP^+/\bfQ^+$.

We have a natural 
identification 
$\phi_{\bfP,\bfQ}:T_\bfQ\is T_\bfQ$
compatible with the natural $\rW_\bfP$-action
where $\rW_\bfP$-acts on $T_\bfQ$ vis the natural embedding 
$\rW_\bfP\to\rW_\bfQ$.

For any $\bfP\in\on{Par}$ and a
non-negative integer $n$, we denote by 
$\bfP_{n}\subset\bfP^+$ the n-th congurence subgroup scheme of $\bfP^+$. Note that we have $\bfP_{0}=\bfP_{}^+$.

There is a bijection between objects in $\on{Par}$
with proper subsets $J\subsetneq\widetilde\Delta$.
For any such $J$ we write $\bfP_J$
for the corresponding parahoric subgroup.
We have $\bfP_{\emptyset}=\bfI$
and $\bfP_{\Delta}=L^+G$.
We write $\bfP^+_J, \rW_J, L_J,$ etc, for the corresponding pro-unipotent radical,  Weyl group, reductive quotient, etc.

\subsection{Grothendieck groups}
For every triangulated category $C$, we denote by
$K_0(C)$ the Grothendieck group of $C$ tensored over $\overline\bbQ_\ell$.
For every object $M\in C$ we denote by $\langle M\rangle\in K_0(C)$
the corresponding isomorphism class.
If $C$ is a monoidal category then $K_0(C)$ is a $\overline\bbQ_\ell$-algebra.
Every triangulated functor $f:C\to C'$
induces a map $\langle f\rangle:K_0(C)\to K_0(C')$.

\section{Sheaf theory admissible ind-scheme and ind-stacks}
\subsection{The case of admissible ind-schemes}
We follow the presentation in \cite[Section 1]{BKV} closely.
Let $\on{Sch}_k$ be the category of quasi-compact and quasi-separated scheme $X$ over $k$ and let $\on{Sch}_k^{ft}$ be the subcategory of separated schemes of finite type over $k$.
For any $X\in\on{Sch}_k$,
let $X/\cdot$ be the category whose objects are morphisms $X\to V$
with $V\in\on{Sch}_k^{ft}$. Following \cite[Section 1.2.2]{BKV},
we denote by
$M(X)$  the colimit
$M(X)=\on{colim}^!_{(X\to V)\in (X/\cdot)^{op}} D(V)$
taken with respect to $!$-pullback
and $D(X)$  the colimit
$D(X)=\on{colim}^*_{(X\to V)\in (X/\cdot)^{op}} D(V)$
taken with respect to $*$-pullback.
Here $D(V)=D_c^b(V,\overline\bbQ_\ell)$ bounded derived category 
of constructible $\overline\bbQ_\ell$-sheaves.
The categories $M(X)$ and $D(X)$ can be viewed 
as the categorical analogs of the space of locally constant measures
and locally constant functions respectively.
For every $f:X\to Y\in\on{Sch}_k$ we have 
natural functors 
$f^!:M(Y)\to M(X)$ and $f^*:D(Y)\to D(X)$.

To define other functors we need the notion of 
admissible schemes and 
admissible morphisms (see \cite[Section 1.1]{BKV}).
A morphism $f:X\to Y\in\on{Sch}_k$ is called admissible 
if there exists a projective system 
$\{X_i\}_{i\in I}$ over $Y$ indexed by a filtered partially ordered set $I$ such that 
each $X_I\to Y$ is finitely presented and all transition maps $X_i\to X_j$, $i>j$ are affine 
unipotent, and $X=\lim_{i\in I} X_i$.
An isomorphism $X=\lim_{i\in I} X_i$ is called an admissible presentation of $X$.
A scheme $X\in \on{Sch}_k$ is called admissible if the structure map
$X\to\on{Spec}(k)$ is admissible.

Let $f:X\to Y$ be an admissible morphism between admissible schemes.
Then we have natural functors 
$f_!:M(X)\to M(Y)$ and $f_*:M(X)\to M(Y)$ (see \cite[Lemma 1.2.4]{BKV}).
Furthermore, if $f$ is of finite presentation then we have functors 
$f^*:M(Y)\to M(X)$, $f_*:M(X)\to M(Y)$, 
$f^!:D(Y)\to D(X)$, $f_!:D(X)\to D(Y)$ (see \cite[Section 1.2.6]{BKV}).

We denote by $\mathbb D_X\in M(X)$ and $1_X\in D(X)$
the dualizing complex and constant sheaf of $X$ respectively.

Let $\on{IndSch}_k$ be the category of ind-schemes over $k$.
Let $X\in\on{Sch}_k$ and $Y\in\on{IndSch}_k$.
A morphism $f:X\to Y$ is called admissible (resp. finitely presented)
if there is a presentation $Y=\on{colim} Y_i$
such that $f$ is induced by an admissible (resp. finitely presented) morphism $f:X\to Y_i$.
If $f$ is a finitely presented closed embedding, then we called 
$X$ is a finitely presented closed subscheme of $Y$.

A morphism $f:X\to Y\in\on{IndSch}_k$ is called admissible (resp. finitely presented) if 
for every finitely presented closed subscheme $Z\subset X$ the map
$f|_Z:Z\to Y$ is admissible (resp. finitely presented).
We call such $f$ schematic if for any finitely presented closed subscheme 
$Z\subset Y$ the preimage $f^{-1}(Z)$ is a scheme.
An ind-scheme $X$ is called admissible if 
the structure map $f:X\to\on{Spec}(k)$ is admissible.

For any ind-scheme $X$, we denote by
$M(X)=\on{colim}_* M(Y)$ and $D(X)=\on{colim}_* D(Y)$ 
where $Y$ runs over the set of finitely presented closed subschemes
and the colimit is taken with respect to
$i_*:M(Y)\to M(Y')$ and $i_*:D(Y)\to D(Y')$.

According to \cite[Section 1.3.2]{BKV},
for every schematic morphism $f:X\to Y$ in $\on{IndSch}_k$
we have functors $f^!:M(Y)\to M(X)$ and $f^*:M(Y)\to M(X)$.
If, in addition, $f$ is finitely presented, we have 
$f^*:M(Y)\to M(X)$. For any admissible morphism $f:X\to Y$ in $\on{IndSch}_k$
we have functors $f_!:M(X)\to M(Y)$ and $f_*:M(X)\to M(Y)$.

Let $X\in\on{IndSch}_k$.
For any finitely presented closed subscheme $Y\subset X$, we denote by
$\delta_Y\in M(X)$ the extension by zero of the dualizing complex $\mathbb D_Y\in M(Y)$.

\subsection{The case of admissible ind-stacks}
Let $\on{St}_k$ be the $2$-category of stacks over $k$
and let $\on{Art}_k^{ft}\subset \on{St}_k$  be the full subcategory of 
Artin stacks of finite type over $k$.
Denote by $\on{St}_k'\subset \on{St}_k$ be the full subcategory 
consisting of $X\in\on{St}_k$ which can be represented by a filtered projective limit
$X\is\lim X_i$ where $X_i\in\on{Art}_k^{ft}$ for all $i$.
For any $X\in\on{Art}_k^{ft}$ one can associate its 
bounded derived category of constructible $\overline\bbQ_\ell$-sheaves 
$D(X)=D_c^b(X,\overline\bbQ_\ell)$.
As explained in \cite[Section 1.4]{BKV}, one can generalize all the notions
defined earlier using $\on{Art}_{k}^{ft}$ and $\on{St}'_{k}$ instead of 
$\on{Sch}^{ft}_k$ and $\on{Sch}_k$, but we do not require the transition morphisms are affine.

For example, we have the notion of 
admissible ind-stacks, admissible morphisms between admissible stacks, 
the categories $M(X)$ and $D(X)$,
and (partially defined) functors 
$f^!, f_!, f^*,f_*$ (see \cite[Section 1.4.5 and 1.4.6]{BKV} for details).

\subsection{Hecke categories}
It is shown in \cite[Section 2.2.2]{BKV}
that the loop group $LG$ and the quotients 
$LG/\bfP_n$, $LG/\bfP_n^+$ are admissible ind-schemes.
We denote by $M(LG)$, $D(LG)$, $M(LG/\bfP_n)$, etc, the corresponding 
categories of sheaves.

It is shown in \cite[Section 2.2.4]{BKV}
that the multiplication map
$m:LG\times LG\to LG$ is admissible and hence we have a functor 
$m_!:M(LG\times LG)\to M(LG)$
and denote by $*$ the convolution 
$\mM*\mM'=m_!(\mM\boxtimes\mM')$.
The convolution $*$ equips $M(LG)$ with a structure of a monoidal category 
(without unit), to be called the Hecke category.

Each $\bfP\in\Par$ acts on $LG$ be the conjugation action 
and we denote by $\frac{LG}{\bfP}$ the corresponding quotient stack.
According to \cite[Section 2.2.5]{BKV}, the quotient stack $\frac{LG}{\bfP}$
is an admissible ind-stack and hence we can form the categories 
$M(\frac{LG}{\bfP})$ and $D(\frac{LG}{\bfP})$.

Consider the following correspondence 
\beq\label{loop groups correspondence}
\frac{LG}{\bfP}\times\frac{LG}{\bfP}\stackrel{\pi}\longleftarrow\frac{LG\times LG}{\bfP}\stackrel{m}\longrightarrow \frac{LG}{\bfP}
\eeq
where $\pi$ is the quotient map and $m$ is multiplication map.
For any $\mM,\mM'\in M(\frac{LG}{\bfP})$ we denote by
$\mM*\mM'=m_!\pi^!(\mM_1\boxtimes\mM_2)$.
The convolution $*$ equips $M(\frac{LG}{\bfP})$ with a structure of monoidal category 
(without unit).

Note that the adjoint action of $\bfP$ on $LG$ descends to an action on
$LG/\bfP^+$ and the quotient 
$\frac{LG/\bfP^+}{\bfP}=\frac{\bfP^+\backslash LG/\bfP^+}{L_\bfP}$ is again a admissible ind-stack.
Denote by 
$M(\frac{LG/\bfP^+}{\bfP})=M(\frac{\bfP^+\backslash LG/\bfP^+}{L_\bfP})$ and $D(\frac{LG/\bfP^+}{\bfP})=D(\frac{\bfP^+\backslash LG/\bfP^+}{L_\bfP})$ the corresponding category of sheaves.
The correspondence in~\eqref{loop groups correspondence} descends to
 the following correspondence 
\[\frac{\bfP^+\backslash LG/\bfP^+}{L_\bfP}\times\frac{\bfP^+\backslash LG/\bfP^+}{L_\bfP}\stackrel{\bar\pi}\longleftarrow\frac{\bfP^+\backslash LG\times^{\bfP^+} LG/\bfP^+}{L_\bfP}\stackrel{\bar m}\longrightarrow \frac{\bfP^+\backslash LG/\bfP^+}{L_\bfP}\]
and we for any $\mM,\mM'\in M(\frac{LG/\bfP^+}{\bfP})=M(\frac{\bfP^+\backslash LG/\bfP^+}{L_\bfP})$ we denote by
$\mM*\mM'=\bar m_!\bar \pi^!(\mM\boxtimes\mM')$.
The convolution $*$ equips $M(\frac{LG/\bfP^+}{\bfP})$ with a structure of monoidal category 
with monoidal unit $\delta_{\frac{\bfP^+/\bfP^+}{\bfP}}$.

We will call the monoidal category 
$(M(\frac{LG/\bfP^+}{\bfP}),*)$ the parahoric Hecke category.

 \section{Stabilization theorem}\label{Bernstein and admissible}
 
\subsection{Induction and restriction functors}
Recall the construction of the 
parabolic induction and 
restriction functors for reductive groups.
Consider the following correspondence
\[\xymatrix{\frac{L}{L}&\frac{P}{P}\ar[r]^{p}\ar[l]_{q}&\frac{G}{G}}\]
The parabolic restriction functors is defined as 
 \[\on{Res}_{L\subset P}^G=q_!p^*:D(\frac{G}{G})\to D(\frac{L}{L})\]
It admits a right adjoint, called the parabolic induction functor, given by
\[\on{Ind}_{L\subset P}^G:=p_*q^!:D(\frac{L}{L})\to D(\frac{G}{G}).\]

We will also consider induction and restriction functors for the non-equivariant 
derived categories. Namely, consider the 
correspondence 
$L \stackrel{\underline q}\la P\stackrel{\underline p}\ra G$ and define the restriction functor as 
\beq\label{res non-equi}
\underline{\on{Res}}_{L\subset P}^G=(\underline q)_!(\underline p)^*:D(G)\to D(L).
\eeq
For the induction functor, consider the 
Grothendieck-Springer alterration 
\[\tilde p:\widetilde G=G\times^B B\ra G,\ \ \ (g,b)\to gbg^{-1}\]
One has natural map 
\[\tilde q:\widetilde G\to T,\ \ \ \ (g,b)\to b\on{\ mod} [B,B]\]
and we define 
\beq\label{Ind non-equi}
\underline{\on{Ind}}_{T\subset B}^G=\tilde p_!\tilde q^*[2\dim U](2\dim U): D(T)\to D(G).
\eeq
We have the following basic compatibility:
\begin{lemma}\label{basic compatibility}
We have
(1) $\pi_T^!\circ\on{Res}_{L\subset P}^G\is \underline{\on{Res}}_{L\subset P}^G\circ\pi_G^![-2\dim U_P](-\dim U_P)
$
and (2)
$\pi_G^!\circ\on{Ind}_{T\subset B}^G\is \underline{\on{Ind}}_{T\subset B}^G\circ\pi_T^![2\dim U](\dim U)$.
\end{lemma}

For any $\mF\in D(T/\rW)$
we will write 
$\mF_{T}\in D(\frac{T}{T})$ (resp. $\underline\mF\in D(T)$)
the $!$-pullback of $\mF$ along the
projection map 
$\frac{T}{T}\is T\times\mathbf B(T)\to T/\rW$
(resp. the map $T\to T/\rW$).
For any $P\in\on{par}$ with Levi subgroup $L$,
it is shown in \cite{C1} that 
the induction $\on{Ind}_{T\subset B_L}^{L}(\mF_{T})$ (resp. 
$\underline\Ind_{T\subset B_L}^{L}(\underline\mF)$)
carries a natural 
$\rW_L$-action and we denote by 
$
\Ind_{T\subset B_L}^{L}(\mF_{T})^{\rW_L}\in D(\frac{L}{L})\ \ 
(\text{resp.}\ \ \underline\Ind_{T\subset B_L}^{L}(\underline\mF)^{\rW_L}\in D(L))
$
the $\rW_L$-invariant summand.

\begin{lemma}\label{W summand}
We have natural isomorphisms $\on{Res}_{L\subset P}^G\on{Ind}_{T\subset B}^{G}(\mF_{T})^{\rW}\is\Ind_{T\subset B_L}^{L}(\mF_{T})^{\rW_L}$
and
 $\underline\Res_{L\subset P}^G\underline\Ind_{T\subset B}^{G}(\underline\mF_{T})^{\rW}\is\underline\Ind_{T\subset B_L}^{L}(\underline\mF_{T})^{\rW_L}
$.
\end{lemma}
\begin{proof}
We give a proof of the first isomorphism.
The second one follows from Lemma \ref{basic compatibility}.
Note that 
the natural adjunction map
$\on{Res}_{L\subset P}^{G}\on{Ind}_{T\subset B}^{G}(-)\to \on{Ind}_{T \subset B_L}^{L}(-)$
gives rise to a natural map
\beq\label{arrow}
\on{Res}_{L\subset P}^{G}\on{Ind}_{T\subset B}^{G}(\mF_{T})^{\rW}
\stackrel{f_1}\to\on{Ind}_{T\subset B_L}^{L}(\mF_T)
\stackrel{f_2}\to \on{Ind}_{T\subset B_L}^{L}(\mF_T)^{\rW_L}
\eeq
where $f_1$ is the induced the adjunction map above and 
$f_2$ is the projection to the $\rW_L$-invariant summand.
In the case when $\mF\in \on{Perv}(T/\rW)$ is a 
perverse sheaf
it follows immediately from
\cite[Proposition 3.4]{C1} that the above map is an isomorphism.

For general $\mF\in D(T/\rW)$,
consider 
the  distinguished triangle 
\[\mF'_{}={^p\tau}_{\leq b-1}(\mF)\to\mF_{}\to\mF''_{}={^p\mathscr H}^b(\mF)[-b]\to\mF'_{}[1]\]
in $D(T/\rW)$,
where $b$ is the largest number such that 
${^p\mathscr H}^b(\mF)\neq 0$. The distinguished triangle above gives rise to the following commutative diagram
\[\xymatrix{\on{Res}_{L\subset P}^{G}\on{Ind}_{T\subset B}^{G}(\mF'_{T})^{\rW}\ar[r]\ar[d]^{h'}&
\on{Res}_{L\subset P}^{G}\on{Ind}_{T\subset B}^{G}(\mF_{T})^{\rW}\ar[r]\ar[d]^h
&\on{Res}_{L\subset P}^{G}\on{Ind}_{T\subset B}^{G}(\mF''_{T})^{\rW}\ar[r]\ar[d]^{h''}&\\
\on{Ind}_{T\subset B_L}^{L}(\mF'_T)^{\rW_L}\ar[r]&\on{Ind}_{T\subset B_L}^{L}(\mF_T)^{\rW_L}\ar[r]&\on{Ind}_{T\subset B_L}^{L}(\mF''_T)^{\rW_L}\ar[r]&}\]
where the vertical arrows are the natural maps in~\eqref{arrow}.
By induction, the maps $h'$ and $h''$ are isomorphisms and it implies that 
$h$ is an isomorphism.

\end{proof}

 \subsection{Harish-Chandra transforms}\label{HC fd}
 We first recall the construction of 
  Harish-Chandra transforms
 for reductive groups and also their relationship with 
 parabolic induction and restriction functors.
 
 For a pair $P\subset Q\in\on{par}$ with 
 unipotent radical $U_P$ and $U_Q$
and Levi quotients $L$ and $M$, one can form the 
 following  horocycle correspondence
 \beq\label{horo}
\frac{G/U_P}{P}\stackrel{h_{P,Q}}\longleftarrow\frac{G/U_Q}{P}\stackrel{c_{P,Q}}\longrightarrow\frac{G/U_Q}{Q}
\eeq
 where $h_{P,Q}$ and $c_{P,Q}$ are the natural maps.
 The  Harish-Chandra transfrom is defined as
\[\on{HC}_{P,Q}=(h_{P,Q})_!(c_{P,Q})^!\is(h_{P,Q})_!(c_{P,Q})^*[2\dim Q/P](\dim Q/P):D(\frac{G/U_Q}{Q})\to
D(\frac{G/U_P}{P}).\]
Since $h_{P,Q}$ is smooth and $c_{P,Q}$ is smooth and proper, the
Harish-Chandra transfrom admits a right adjoint given by
\[\on{CH}_{P,Q}=(c_{P,Q})_!(h_{P,Q})^*:D(\frac{G/U_P}{P})\to
D(\frac{G/U_Q}{Q}).\]

We have the following basic properties of 
Harish-Chandra transforms and their compatibility with 
parabolic restriction functors.

\begin{lemma}\label{compatibility}
(1) The functor $\on{HC}_{P,Q}$ is monoidal. 

(2)  For any $P\subset Q\subset Q'$ parabolic subgroups,
we have natural isomorphisms of functors 
\[\on{CH}_{Q,Q'}\circ\on{CH}_{P,Q}\is\on{CH}_{P,Q'}\ \ \ \ \on{HC}_{P,Q}\circ\on{HC}_{Q,Q'}\is\on{HC}_{P,Q'}\]

(3)
Denote by $f_P:\frac{L}{P}\to\frac{L}{L}$ the natural projection map
and $i_P:\frac{L}{P}\to\frac{G/U}{P}$ the closed embedding.
There are canonical isomorphisms of functors
\[(f_P)_!(i_P)^*\on{HC}_{P,G}\is \on{Res}_{L\subset P}^G[2\dim U_P](\dim U_P)\]
\[\on{CH}_{P,G}(i_P)_!f_P^!\is \on{Ind}_{L\subset P}^G(\mF)[-2\dim U_P](-\dim U_P)\]

(4)
For any $P\subset Q\in\on{par}$, $\mF_P\in D(\frac{G/U_P}{P})$
and $\mF_Q\in D(\frac{G/U_Q}{Q})$, we have 
\[\mF_Q*\on{CH}_{P,Q}(\mF_P)\is\on{CH}_{P,Q}(\on{HC}_{P,Q}(\mF_Q)*\mF_P)\]

\end{lemma}
\begin{proof}
Part (1), (2), and (3) are standard facts.
For part (4),
we give a proof when $P\subset G\in\on{par}$
(since this is the case that will be used later).
 The proof in the general case is 
similar.
Consider the following diagram
\[\xymatrix{&\frac{G}{P}\ar[ld]_c\ar[r]^h&\frac{G/U_P}{P}&&\\
\frac{G}{G}\ar[d]&\frac{G\times G}{P}\ar[u]\ar[l]\ar[d]\ar[r]&\frac{G\times (G/U_P)}{P}\ar[r]_b\ar[u]^a&\frac{G}{G}\times\frac{G/U_P}{P}\ar[d]\\
\frac{G}{G}&\frac{G\times G}{P}\ar[l]\ar[d]\ar[r]&\frac{G}{G}\times\frac{G}{P}\ar[d]^{\id\times c}\ar[r]^{\id\times h}&\frac{G}{G}\times\frac{G/U_P}{P}\\
&\frac{G\times G}{G}\ar[r]^q\ar[lu]^m&\frac{G}{G}\times\frac{G}{G}&}\]
where $a$ and $m$ are the multiplication maps and other map 
are the natural  quotients maps.
The top and bottom square diagrams are Cartesian and using base-change theorems one finds that 
\[\mF_G*\on{CH}_{P,G}(\mF_P)\is m_!q^!(\id\times c)_!(\id\times h)^*(\mF_G\times\mF_P)\is
c_!h^*a_! b^!(\mF_G\boxtimes\mF_P)\is\on{CH}_{P,G}(a_! b^!(\mF_G\boxtimes\mF_P)).\]
Consider the diagram
\[\xymatrix{\frac{G/U_P}{P}&\frac{G\times G/U_P}{P}\ar[r]\ar[l]_a\ar[d]\ar@/^0.6cm/[rr]^b&\frac{G}{P}\times\frac{G/U_P}{P}\ar[r]_{c\times\id}\ar[d]^{h\times\id}&\frac{G}{G}\times\frac{G/U_P}{P}\\
&\frac{U_P\backslash G\times^{U_P}G/U_P}{L_P}\ar[r]^l\ar[lu]^m&\frac{G/U_P}{P}\times
\frac{G/U_P}{P}}.\]
Again by base change theorems, one finds that 
\[a_! b^!(\mF_G\boxtimes\mF_P)\is m_!l^!(h\times\id)_!(c\times\id)!(\mF_G\times\mF_Q)\is\on{HC}_{P,Q}(\mF_G)*\mF_P\]
All together, we obtain the desired isomorphism
\[\mF_G*\on{CH}_{P,G}(\mF_P)\is\on{CH}_{P,G}(a_! b^!(\mF_G\boxtimes\mF_P))\is\on{CH}_{P,G}(\on{HC}_{P,Q}(\mF_G)*\mF_P).\]

\end{proof}

Consider the parabolic 
Springer map 
\[s_{P,Q}:\tilde\calN_{P,Q}=\{(x,gP)\in M\times Q/P|x\in gU_Pg^{-1}\}/Q\to\frac{M}{Q}\ \ \ s_{P,Q}(x,gP)=x\]
and the corresponding 
parabolic Springer sheaf
\[\mathcal{S}_{P,Q}=(s_{P,Q})_!\overline{\bbQ}_\ell[2\dim U_P\cap L_Q]\in D(\frac{M}{Q})\]
Denote by
$i_Q:\frac{M}{Q}\is\frac{Q/U_Q}{Q}\to\frac{G/U_Q}{Q}$
the natural inclusion map.
We have the following important property of 
the Harish-Chandra transform.

\begin{proposition}
For any $\mF\in D(\frac{G/U_Q}{Q})$, we have 
$\on{CH}_{P,Q}\circ\on{HC}_{P,Q}(\mF)\is((i_{Q})_!\calS_{P,Q})*\mF$.
In particular, the identify functor is a direct summand of $\on{CH}_{P,Q}\circ\on{HC}_{P,Q}$ and $\on{HC}_{P,Q}$ is conservative.
\end{proposition}
\begin{proof}
This is proved in \cite{G} in the case when $P=B, Q=G$, and the general case in \cite[Proposition 4.1]{BT}.
\end{proof}

 \begin{corollary}\label{B implies P}
 Let $\mF\in D(\frac{G}{G})$. Assume 
 $\on{HC}_{B,G}(\mF)$ is supported on 
 $\frac{B/U}{B}\subset\frac{G/U}{B}$.
 \begin{enumerate}
\item
Then for any standard parabolic subgroup $P$, we have 
\beq\label{formula HC=Res}
\on{HC}_{P,G}(\mF)\is (i_P)_*(f_P)^!\on{Res}_{L\subset P}^G(\mF)[2\dim U_P](\dim U_P).
\eeq
 In particular, $\on{HC}_{P,G}(\mF)$ is supported on 
$\frac{P/U_P}{P}\subset \frac{G/U_P}{P}$
 \item
 For any $\mF_T\in D(\frac{T}{T})$, we have 
 \[\mF*\Ind_{T\subset B}^G(\mF_T)\is\Ind_{T\subset B}^G(\Res_{T\subset B}^G(\mF)*\mF_T)[2\dim U](\dim U).\]
\end{enumerate}
 \end{corollary}
 \begin{proof}
 Proof of (1).
We have the following commutative diagram
\beq\label{supp}
\xymatrix{\frac{B/U}{B}\ar[d]^{i_B}&\frac{B/U_P}{B}\ar[r]\ar[l]\ar[d]&\frac{P/U_P}{P}\ar[d]\\
\frac{G/U}{B}&\frac{G/U_P}{B}\ar[r]^{c_{B,P}}\ar[l]_{h_{B,P}}&\frac{G/U_P}{P}}\eeq
where the vertical maps are the natural inclusion 
and the left square is Cartesian. 
Write $\mF'=\on{HC}_{P,G}(\mF)$.
Since 
\[\on{HC}_{B,P}(\mF')\is
\on{HC}_{B,P}\circ\on{HC}_{P,G}(\mF)\is\on{HC}_{B,G}(\mF) \] is supported on $\frac{B/U}{B}$, 
we have 
\[\on{HC}_{B,P}(\mF')\is
(i_B)_!(i_B)^*\on{HC}_{B,P}(\mF')\]
and the  base change formula implies that 
\[\on{CH}_{B,P}\circ\on{HC}_{B,P}(\mF')\is\on{CH}_{B,P}\circ(i_B)_!(i_B)^*\on{HC}_{B,P}(\mF')\is (c_{B,P})_!(h_{B,P})^*
(i_B)_!(i_B)^*\on{HC}_{B,P}(\mF)\] 
is supported $\frac{P/U_P}{P}$.
Since $\mF'$ is a summand of $\on{CH}_{B,P}\circ\on{HC}_{B,P}(\mF')$,
it implies that $\mF'$ is also supported on 
$\frac{P/U_P}{P}$.
The desired isomorphism~\eqref{formula HC=Res} follow from part (3) of Lemma \ref{compatibility}.

Part (2) follows from part (1) and Corollary \ref{compatibility}.
 \end{proof}

%%%%%
\quash{
We have the following commutative diagram
\beq\label{supp}
\xymatrix{\frac{P/U_P}{P}\ar[d]^{i_P}&\frac{P/U_Q}{P}\ar[r]\ar[l]\ar[d]&\frac{Q/U_Q}{Q}\ar[d]\\
\frac{G/U_P}{P}&\frac{G/U_Q}{P}\ar[r]^{c_{P,Q}}\ar[l]_{h_{P,Q}}&\frac{G/U_Q}{Q}}\eeq}
%%%%%

\subsection{Harish-Chandra transforms for loop groups}\label{affine HC}

For any $\bfP\subset\bfQ\in\on{Par}$,
we can form the following  horocycle correspondence for loop groups:
\beq\label{affine horo}
\frac{LG/\bfP^+}{\bfP}\stackrel{h_{\bfP,\bfQ}}\longleftarrow\frac{LG/\bfQ^+}{\bfP}\stackrel{c_{\bfP,\bfQ}}\longrightarrow\frac{LG/\bfQ^+}{\bfQ}
\eeq
where the maps $h_{\bfP,\bfQ}$ and $c_{\bfP,\bfQ}$ are the natural projection maps.
The  Harish-Chandra transfrom is defined as
\[\on{HC}_{\bfP,\bfQ}:=(h_{\bfP,\bfQ})_!(c_{\bfP,\bfQ})^!:M(\frac{LG/\bfQ^+}{\bfQ})\to
M(\frac{LG/\bfP^+}{\bfP}).\]
Since $h_{\bfP,\bfQ}$ is smooth and $c_{\bfP,\bfQ}$ is smooth and proper, the
Harish-Chandra transfrom admits a right adjoint given by
\[\on{CH}_{\bfP,\bfQ}:=(c_{\bfP,\bfQ})_!(h_{\bfP,\bfQ})^*:M(\frac{LG/\bfP^+}{\bfP})\to
M(\frac{LG/\bfQ^+}{\bfQ}).\]
Consider the projection maps
\[
\pi_\bfP:LG\to\frac{LG/\bfP^+}{\bfP}\ \ \ \ \ \ 
 h_\bfP:\frac{LG}{\bfP}\to\frac{LG/\bfP^+}{\bfP}
 \ \ \ \ \ \
 \pi_{\bfP,\bfQ}:\frac{LG}{\bfP}\to\frac{LG}{\bfQ}
.\]

%Consider the following projection maps
%\beq\label{projection}
%\pi_\bfP:LG\to\frac{LG/\bfP^+}{\bfP}\ \ \ \ \ \
%h_{\bfP}:\frac{LG}{\bfP}\to\frac{LG/\bfP^+}{\bfP}.
%\eeq

\begin{proposition}\label{aHC}
(1)
For any $\mathbf{P}\subset\mathbf{Q}\subset\mathbf{Q}'\in\Par$
We have natural isomorphisms of functors 
$\on{CH}_{\mathbf{Q},\mathbf{Q}'}\circ\on{CH}_{\mathbf{P},\mathbf{Q}}\is\on{CH}_{\mathbf{P},\mathbf{Q}'}$
and $\on{HC}_{\mathbf{P},\mathbf{Q}}\circ\on{HC}_{\mathbf{Q},\mathbf{Q}'}\is\on{HC}_{\mathbf{P},\mathbf{Q}'}$

(2) 
The functor $\on{HC}_{\mathbf{P},\mathbf{Q}}$ is monoidal, that is,
there is a canonical isomorphism 
\[\on{HC}_{\mathbf{P},\mathbf{Q}}(\mM*\mM')\is\on{HC}_{\mathbf{P},\mathbf{Q}}(\mM)*\on{HC}_{\mathbf{P},\mathbf{Q}}(\mM')\]
satisfying the natural compatibility conditions. 

%(3)
%For any $\mM\in M(\frac{LG/\mathbf Q^+}{\mathbf Q})$, we have 
%$\on{CH}_{\mathbf P,\mathbf Q}\circ\on{HC}_{\mathbf P,\mathbf Q}(\mM)\is((i_{\mathbf Q})_*\calS_{P,Q})*\mM$.
%In particular, the identify functor is a direct summand of $\on{CH}_{\mathbf P,\mathbf Q}\circ\on{HC}_{\mathbf P,\mathbf Q}$ and $\on{HC}_{\mathbf P,\mathbf Q}$ is conservative.

(3)
For any $\mM\in M(\frac{LG/\mathbf{Q}^+}{\mathbf{Q}})$, we have 
\beq\label{convolution 1}
\pi_\mathbf{P}^!(\on{HC}_{\mathbf{P},\mathbf{Q}}(\mM))
\is
\pi_{\mathbf{Q}}^!(\mM)*\delta_{\mathbf{P}^+}
\eeq
\beq\label{convolution 2}
h_\mathbf P^!\on{HC}_{\mathbf P,\mathbf Q}(\mM)\is \pi_{\mathbf P,\mathbf Q}^!h_\bfQ^!(\mM)*\delta_{\frac{\mathbf P^+}{\mathbf P}}.
\eeq

\end{proposition}
\begin{proof}
The proof of part (1) and (2) are similar to the finite dimensional case and we omit the details.
For part (3), we give a proof of the first isomorphism. The proof of the second one is  similar.
Consider the following diagram
\[\xymatrix{&LG\times\bfP^+\ar[d]^{}\ar[d]^{q_1}\ar[rd]^{\pr}\ar[ld]_{a}&\\
LG\ar[d]^{\pi_\bfP}&LG\times^{\bfQ^+}\bfP^+\ar[l]_{\bar a}\ar[d]\ar[rd]^{q_2}&LG\ar[d]^{\pi_\bfQ}\\
\frac{LG/\bfP^+}{\bfP}&\frac{LG/\bfQ^+}{\bfP}\ar[l]_{h_{\bfP,\bfQ}}\ar[r]^{c_{\bfP,\bfQ}}&\frac{LG/\bfQ^+}{\bfQ}}\]
where $a$ and $\bar a$ are multiplication maps and $q_1$ and $q_2$ are the quotient maps.
Since the left corner square is Cartesian and $\pi_\bfP$ is formally smooth, the base change theorems
imply that 
$\pi_\mathbf{P}^!\on{HC}_{\mathbf{P},\mathbf{Q}}(\mM)\is \pi_\mathbf{P}^! (h_{\bfP,\bfQ})_! c_{\bfP,\bfQ}^!(\mM)\is \bar a_!q_2^!(\mM)$. Note that $(q_1)_!q_1^!\is\id$ and hence $\bar a_!q_2^!(\mM)\is\bar a_!(q_1)_!q_1^!q_2^!(\mM)\is a_!\pr^!\pi_\bfQ^!(\mM)\is\pi_\bfQ^!(\mM)*\delta_{\bfP^+}$.
The desired isomorphism follows.
\end{proof}

%\begin{remark}
%Unlike the Harish-Chandra transform, the adjoint functors 
%$\on{CH}_{P,Q}$ or
%$\on{CH}_{\bfP,\bfQ}$ are not monoidal.
%\end{remark}

 \subsection{The algebra $\mathfrak A(LG)$}\label{admissible collections}
 Consider the following space 
\[\frakA(LG):=\on{lim}_{\bfP\in\Par} K_0(M(\frac{LG/\bfP^+}{\bfP}))\]
where the limit is taken with respect to the  Harish-Chandra transforms
\[\langle\on{HC}_{\bfP,\bfQ}\rangle:K_0(M(\frac{LG/\bfQ^+}{\bfQ}))\to
K_0(M(\frac{LG/\bfP^+}{\bfP}))\]
Objects of $\frakA(LG)$
are collections $\mM=\{\langle\mM_\bfP\rangle\}_{\bfP\in\Par}$ 
such that $\langle\on{HC}_{\bfP,\bfQ}\rangle(\langle\mM_\bfQ\rangle)=\langle\mM_\bfP\rangle$.

\begin{example}\label{unit}
For any $\on{P}\in\Par$ let $\delta_{\frac{\bfP^+/\bfP^+}{\bfP}}$
the dualizing complex on $\frac{\bfP^+/\bfP^+}{\bfP}$. 
Then using Lemma \ref{aHC} (3) and the fact that 
$\delta_{\bfP^+}*\delta_{\bfQ^+}\is\delta_{\bfP^+}$
for $\bfP\subset\bfQ\in\Par$, we see that  
$\delta:=\{\langle\delta_{\frac{\bfP^+/\bfP^+}{\bfP}}\rangle\}_{\bfP\in\Par}$
is in $\frakA(LG)$.
\end{example}

The monoidal structures on the affine Harish-Chandra transforms 
in Lemma \ref{aHC}, implies

\begin{lemma}
There is a
natural unital algebra structure on $\frakA(LG)$
where the multiplication is given by
$\mM*\mM':=\{\langle\mM_\mathbf{P}*\mM_{\mathbf{P}'}\rangle\}_{\mathbf P\in\Par}$
and the unit is $\delta$.
\end{lemma}

\subsection{The subalgebra $\frakA(LG)_1$}
The closed embedding 
$\frac{\bfP/\bfP^+}{\bfP}\to \frac{LG/\bfP^+}{\bfP}$
induces  an injective map
$K_0(M(\frac{\bfP/\bfP^+}{\bfP}))\to K_0(M(\frac{LG/\bfP^+}{\bfP}))$
and we can form the 
the following subspace  
of $\frakA(LG)$:
\beq\label{A(LG)_1}
\frakA(LG)_1=
\{\{\langle\mM_\bfP\rangle\}_{\bfP\in\Par}\in\frakA(LG)|\langle\mM_\bfP\rangle\in K_0(M(\frac{\bfP/\bfP^+}{\bfP}))\}
\eeq
Since $\bfP\subset LG$
is a subgroup, the subspace $K_0(M(\frac{\bfP/\bfP^+}{\bfP}))\subset K_0(M(\frac{LG/\bfP^+}{\bfP}))$ is in fact a subalgebra
and it follows that $\frakA(LG)_1$ is also a subalgebra
$\frakA(LG)$.  

We shall give an alternative description of 
$\frakA(LG)_1$ in terms of sheaves on the Levi quotients.
To this end, let us consider the following subcategory 
of $D(\frac{G}{G})$:
\beq\label{A(G)}
A(G)=\{\mF\in D(\frac{G}{G})| \on{Supp}(\on{HC}_{P,G}(\mF))\subset\frac{P/U_P}{P} \ \text{for all }\  P\in\on{par} \}
\eeq
Note that $A(G)$ is a monoidal subcategory of $D(\frac{G}{G})$.

\begin{lemma}\label{par}
Let $\mF\in A(G)$
and $P\in\on{par}$ with Levi subgroup $L$.
(1)
We have 
$\on{Res}_{L\subset P}^G(\mF)\in A(L)$.
(2) The functor
$\on{Res}_{L\subset P}^G[2\dim U_P](\dim U_P):D(\frac{G}{G})\to D(\frac{L}{L})$
is monoidal.
%For any $\mF'\in D(\frac{G}{G})$, we have   
%$\on{Res}_{L\subset P}^G(\mF*\mF')\is\on{Res}_{L\subset P}^G\mF*\on{Res}_{L\subset P}^G\mF'[2\dim U_P](\dim U_P)$.
\end{lemma}
\begin{proof}
Part (2) follows from the monoidal property of 
Harish-Chandra functor and~\eqref{formula}.
Proof of (1).
Write $B_L=B/U_P\subset P/U_P=L$
for the Borel subgroup of $L$
and $\mF'=\on{Res}_{L\subset P}^G(\mF)$.
By Corollary \ref{B implies P}, it suffices to show that 
$\on{HC}_{B_L,L}(\mF')$
is supported on $\frac{B_L/U_L}{B_L}$.
Consider the following 
cartesian diagrams
\beq\label{supp}
\xymatrix{\frac{L/U_L}{B_L}&\frac{L}{B_L}\ar[r]\ar[l]&\frac{L}{L}\\
\frac{L/U_L}{B}\ar[d]_{u}\ar[u]^v&\frac{L}{B}\ar[r]\ar[l]\ar[d]\ar[u]
&\frac{L}{P}\ar[d]_{i_P}\ar[u]^{f_P}&\frac{P}{P}\ar[l]\ar[r]\ar[d]&\frac{G}{G}\ar[d]\\
\frac{G/U}{B}&\frac{G/U_P}{B}\ar[r]^{c_{B,P}}\ar[l]_{h_{B,P}}&\frac{G/U_P}{P}
&\frac{G}{P}\ar[r]\ar[l]&\frac{G}{G}}\eeq
where the arrows are the natural inclusions and quotient maps.
By applying the base-change formula to the above diagram
, we see that 
\[v^!\on{HC}_{B_L,L}(\mF')\is
u^*\on{HC}_{B,P}\circ (i_P)*f_P^!\mF'
\stackrel{~\eqref{formula HC=Res}}\is u^*\on{HC}_{B,P}\circ\on{HC}_{P,G}(\mF)[-2\dim U_P]\is\]
\[\is u^*\on{HC}_{B,G}(\mF)[-2\dim U_P]\]
Since $\on{HC}_{B,G}(\mF)$ is supported on $\frac{T}{B}$, the above isomorphism implies 
$v^!\on{HC}_{B_L,L}(\mF')$ 
is supported on $\frac{B_L/U_L}{B}$
and hence $\on{HC}_{B_L,L}(\mF')$
is supported on $\frac{B_L/U_L}{B_L}$.
\end{proof}

The lemma above implies that for any $\bfP\subset\bfQ\in\Par$
the  parabolic restriction functors 
restricts to a monoidal functor
\beq\label{res}
\on{Res}_{L_\bfP,B_{\bfP,\bfQ}}^{L_\bfQ}[2\dim U_{\bfP,\bfQ}](\dim U_{\bfP,\bfQ}):A(L_\bfQ)\to A(L_\bfP).
\eeq
Consider the following algebra
\beq
\lim_{\bfP\in\on{Par}}K_0(A(L_\bfP))
\eeq
where the limit is taken with respect to parabolic restriction functors
~\eqref{res}.

For any $\bfP\in\Par$ 
, let $f_\bfP:\frac{\bfP/\bfP^+}{\bfP}\to \frac{L_\bfP}{L_\bfP}$
be the projection map.

 \begin{lemma}\label{1=par}
 There is 
an algebra isomorphism 
$\upsilon:\lim_{\mathbf P\in\on{Par}}K_0(A(L_\mathbf P))
\is\frakA(LG)_1$ given by
\[\upsilon(\{\langle\mF_{L_{\mathbf P}}\rangle\}_{\mathbf P\in\Par})=\{\langle f_\mathbf P^!\rangle(\langle\mF_{L_{\mathbf P}}\rangle)\}_{\mathbf P\in\Par}\]

\end{lemma}
\begin{proof}
We first show that the map $\upsilon$ is well-defined.
There is a commutative diagram
\beq\label{finite horocycle}
\xymatrix{&\frac{\bfP/\bfP^{+}}{\bfP}\ar[r]\ar[d]\ar[ld]_{f_\bfP}&\frac{\bfQ/\bfP^{+}}{\bfP}\ar[d]^{f_{\bfP,\bfQ}}&\frac{\bfQ/\bfQ^+}{\bfP}\ar[r]\ar[d]\ar[l]&\frac{\bfQ/\bfQ^+}{\bfQ}\ar[d]^{f_\bfQ}\\
\frac{L_\bfP}{L_\bfP}&\frac{L_\bfP}{B_{\bfP,\bfQ}}\ar[r]_{i_{B_{\bfP,\bfQ}}}\ar[l]^{f_{B_{\bfP,\bfQ}}}&
 \frac{L_{\bfQ}/U_{\bfP,\bfQ}}{B_{\bfP,\bfQ}}&\frac{L_{\bfQ}}{B_{\bfP,\bfQ}}\ar[r]\ar[l]&\frac{L_{\bfQ}}{L_{\bfQ}}}
 \eeq
 where the upper horizontal arrows are the restriction of 
 the maps $c_{\bfP,\bfQ}$ and $h_{\bfP,\bfQ}$ to the 
 closed subscheme $\frac{\bfQ/\bfQ^+}{\bfP}\subset
 \frac{LG/\bfQ^+}{\bfP}$, 
 the lower horizontal arrows are the 
horocycle correspondence in~\eqref{horo}
associated to the parabolic subgroup 
$B_{\bfP,\bfQ}\subset L_{\bfQ}$.
Note also the all the vertical map are torsors over pro-unipotent groups
and hence $!$-pushforward along those maps are equivalences with inverse 
equivalences given by $!$-pullback.
It follows that 
\[\on{HC}_{\bfP,\bfQ}(f_\bfQ^!\mF_{L_{\bfQ}})\is  f_{\bfP,\bfQ}^!\on{HC}_{B_{\bf P,\bfQ},L_\bfQ}(\mF_{L_\bfQ}))\is\]
\[\stackrel{\eqref{formula}}\is f_{\bfP,\bfQ}^!(i_{B_{\bfP,\bfQ}})_*(f_{B_{\bfP,\bfQ}})^!\on{Res}_{L_\bfP\subset B_{\bfP,\bfQ}}^{L_{\bfQ}}(\mF_{L_\bfQ})[2\dim U_{\bfP,\bfQ}](\dim U_{\bfP,\bfQ})\is f_\bfP^!(\mF_{L_\bfP}).\]
The claim follows.
The inverse  is given by
$\upsilon^{-1}(\{\langle\mF_\bfP\rangle\}_{\bfP\in\Par})=\{\langle(f_\bfP)_!\rangle\langle\mF_{\bfP}\rangle\}_{\bfP\in\Par}$.

\end{proof}

\begin{remark}
At the moment, it is not clear whether the algebra
$\frakA(LG)_1$ contains any  object
other than the unit element $\delta$. 
In Section \ref{SCC}, we will use the 
Lemma \ref{1=par}
and the  results in \cite{C2} to produce many interesting objects in $\frakA(LG)_1$.
\end{remark}

 \subsection{Averaging functors}
 Let $a:LG\times LG\to LG$ be the adjoint action map
$a(x,y)=xyx^{-1}$.
For any $\mM\in M(LG)$ gives rise to a functor
 \[\on{Av}^\mM:M(LG)\to M(LG)\ \ \ \ \on{Av}^\mM(\mF)=a_!(\mM\boxtimes\mF).\]
 If $\mM=1_X$ for some closed subscheme $X\subset LG$
 we will write 
 $\on{Av}^X=\on{Av}^{1_X}$.
 
For any parahoric $\bfP\in\on{Par}$ we have the following maps
\[LG/\bfP\times\frac{LG}{\bfP}\stackrel{q}\longleftarrow LG\times^{\bfP}LG\stackrel{a}\longrightarrow LG\]
where $q$ is the quotient map. Let $\mM\in M(\frac{LG}{\bfP})$
 and let $Y\subset LG/\bfP$ be a locally closed subscheme.
 We define the averaging functor to be 
 \[\on{Av}^Y:M(\frac{LG}{\bfP})\to M(LG)\ \ \ \on{Av}^Y(\mM)=a_!q^!( 1_Y\boxtimes\mM)\]
 where $1_Y\in M(LG/\bfP)=D(LG/\bfP)$ is the constant sheaf 
 on $Y$.
 
 Assume $Y$ is invariant under the action of $\bfQ\in\on{Par}$. 
We can consider the quotient 
\[\bfQ\backslash LG/\bfP\times\frac{LG}{\bfP}\stackrel{ q}\longleftarrow \bfQ\backslash LG\times^{\bfP}LG\stackrel{a}\longrightarrow\frac{LG}{\bfQ}\] 
and we define the following $\bfQ$-equivariant version of the averaging functors:
 \[\on{Av}^Y_{\bfQ}:M(\frac{LG}{\bfP})\to M(\frac{LG}{\bfQ})\ \ \ \on{Av}^Y_{\bfQ}(\mM)=a_! q^!( 1_{\bfQ\backslash Y}\boxtimes\mM)\]
 where $1_{\bfQ\backslash Y}\in M(\bfQ\backslash LG/\bfP)$ is the 
 unique measure such that its $!$-pull-back to $M(LG/\bfP)$ is 
 $1_Y$.

\begin{example}
Assume $\bfP\subset\bfQ$
and 
 $Y=\bfQ/\bfP\subset LG/\bfP$, which is 
 $\bfQ$-invariant. Let $\pi_{\bfP,\bfQ}:\frac{LG}{\bfP}\to\frac{LG}{\bfQ}$
 be the natural projection map.
It is straightforward to check that  
 \[\on{Av}_{\bfQ}^{\bfQ/\bfP}\is(\pi_{\bfP,\bfQ})_![2\dim Y]:M(\frac{LG}{\bfP})\to M(\frac{LG}{\bfQ})\]
 \end{example}

 \subsection{Stabilization theorem for objects in $\frakA(LG)$}
 Let $\Upsilon$ be the partially ordered set 
of non-empty $\bfI$-invariant subscheme $\mathrm Y\subset LG/\bfI$.
For any $\mathrm Y\in\Upsilon$ we write $\widetilde{\mathrm Y}\subset LG$ its preimage in
$LG$. For any $\mathrm Y\in\Upsilon$ and $\bfP\in\Par$ we write 
$\mathrm Y_\bfP=\widetilde {\mathrm Y}\bfP/\bfP\subset LG/\bfP$.

Fix an object $\mM=\{\langle\mM_\bfP\rangle\}_{\bfP\in\Par}\in\frakA(LG)$.
For any $\mathrm Y\in\Upsilon$, we set
\beq\label{A^Y_M}
\langle A^\mathrm Y_\mM\rangle:=\sum_{\bfP\in\Par}(-1)^{r(G)-r(\bfP)}\langle\on{Av}^{\mathrm Y_\bfP}(h_\bfP^!\mM_\bfP)\rangle\in K_0(M(LG)).
\eeq

We shall prove the following stabilization theorem.
\begin{thm}\label{stabilization}
(1)
For every $\mF\in M(LG)$, the system $\{\langle A^\mathrm Y_\mM\rangle*\langle\mF\rangle\}_{\mathrm Y\in\Upsilon}$
stabilizes.

(2) For every $\mathbf{P}\in\Par$ and $\mathrm Y\in\Upsilon$, we have 
$\langle A^\mathrm Y_\mM\rangle*\langle\delta_{\mathbf{P}^+}\rangle=\langle\pi^!_{\mathbf{P}}\mM_\mathbf{P}\rangle$

\end{thm}
The theorem above generalizes
\cite[Theorem 4.1.5]{BKV} in the case when $\mM=\delta\in\frakA(LG)$ is the unit (Example \ref{unit}). It tuns out that the argument given in \emph{loc. cit.}
is robust enough to treat the case of objects in $\frakA(LG)$.
The key step is the following generalization of 
\cite[Lemma 4.2.3]{BKV}.
For any $w\in\widetilde\rW$, let 
$\mathrm Y_w=\bfI w\bfI/\bfI\subset LG/\bfI$
and $J_w=\{\alpha\in\widetilde\Delta|w(\alpha)>0\}$.
\begin{lemma}\label{key lemma}
Let $w\in\widetilde\rW$, $\alpha\in\widetilde\Delta$, $\mathbf Q\in\Par$,
$n\in\mathbb N$, and $J\subset J_{w}\setminus\alpha$
be such that $U_{w(\alpha)}\in\mathbf Q_n^+$.
Write $J'=J\cup\alpha$.
We have 
\[\langle\on{Av}^{(Y_w)_{\mathbf P_{J'}}}(h_{\mathbf P_{J'}}^!\mM_{\mathbf P_{J'}})\rangle*\langle\delta_{\mathbf Q_n^+}\rangle=\langle\on{Av}^{(Y_w)_{\mathbf P_{J}}}(h_{\mathbf P_{J}}^!\mM_{\mathbf P_{J}})\rangle*\langle\delta_{\mathbf Q_n^+}\rangle\]
\end{lemma}
\begin{proof}
Let $\beta_1,...,\beta_{l(w)}$ be all positive affine roots such 
that $w(\beta_i)$ is an affine negative root  and set $\bfI_w=\prod_{i=1}^{l(w)} U_{\beta_i}\subset\bfI^+$.
Pick a representative $\underline w\in LG$ of $w\in\widetilde\rW$
and consider the closed subscheme $\bfI_w\underline w\subset LG$. 
Then the same proof of the isomorphism (4.3) of \cite{BKV}, using 
\cite[Lemma 2.3.2 (a)]{BKV}, shows that we have an equality
\[\langle\on{Av}^{(Y_w)_{\mathbf P_{J}}}\rangle(\langle h_{\mathbf P_{J}}^!\mM_{\mathbf P_{J}}\rangle)*\langle\delta_{\mathbf Q_n^+}\rangle= 
\langle\on{Av}^{\bfI_w\underline w}\rangle(\langle\pi_{\mathbf P_J}^!\mM_{\mathbf P_J}\rangle*\langle\delta_{\underline w^{-1}\mathbf Q_n^+\underline w} \rangle)\]
and similarly for $J'$.
Thus it suffices to show that there an equality
\beq\label{key iso}
\langle\pi_{\mathbf P_{J'}}^!\mM_{\mathbf P_{J'}}\rangle*\langle\delta_{\underline w^{-1}\mathbf Q_n^+\underline w}\rangle=\langle\pi_{\mathbf P_J}^!\mM_{\mathbf P_J}\rangle*\langle\delta_{\underline w^{-1}\mathbf Q_n^+\underline w}\rangle
\eeq
By \cite[Lemma 4.2.3]{BKV}, the natural map 
$\delta_{\bfP_{J'}^+}\to\delta_{\bfP_J^+}$ induces an equality 
\beq\label{step 1}
\langle\delta_{\bfP_{J'}^+}\rangle*\langle\delta_{\underline w^{-1}\mathbf Q_n^+\underline w}\rangle=\langle
\delta_{\bfP_J^+}\rangle*\langle\delta_{\underline w^{-1}\mathbf Q_n^+\underline w}\rangle.
\eeq
On the other hand, we have
 \beq\label{step 2}
 \langle\pi_{\mathbf P_{J'}}^!\mM_{\mathbf P_{J'}}\rangle=\langle\pi_{\mathbf P_{J'}}^!\mM_{\mathbf P_{J'}}\rangle*\langle\delta_{\bfP_{J'}^+}\rangle
 \eeq
and, by Proposition \ref{aHC},  
\beq\label{step 3}
\langle\pi_{\mathbf P_{J'}}^!\mM_{\mathbf P_{J'}}\rangle*\langle\delta_{\bfP_J^+}\rangle\stackrel{~\eqref{convolution 1}}=
\langle\pi_{\bfP_J}^!\rangle\langle\on{HC}_{\bfP_J,\bfP_{J'}}(\mM_{\bfP_{J'}})\rangle=
\langle\pi_{\bfP_J}^!
\rangle\langle\mM_{\bfP_J}\rangle,
\eeq
where the last equality follows from the definition of $\frakA(LG)$.
All together, we obtain the desired equality
\[\langle\pi_{\mathbf P_{J'}}^!\mM_{\mathbf P_{J'}}\rangle*\langle\delta_{\underline w^{-1}\mathbf Q_n^+\underline w}\rangle\stackrel{~\eqref{step 2}}=\langle
\pi_{\mathbf P_{J'}}^!\mM_{\mathbf P_{J'}}\rangle*\langle\delta_{\bfP_{J'}^+}\rangle*\langle\delta_{\underline w^{-1}\mathbf Q_n^+\underline w}\rangle\stackrel{~\eqref{step 1}}=\langle
\pi_{\mathbf P_{J'}}^!\mM_{\mathbf P_{J'}}\rangle*\langle\delta_{\bfP_{J}^+}\rangle*\langle\delta_{\underline w^{-1}\mathbf Q_n^+\underline w}\rangle= \]
\[\stackrel{~\eqref{step 3}}=
\langle\pi_{\mathbf P_{J}}^!\mM_{\mathbf P_{J}}\rangle*\langle\delta_{\underline w^{-1}\mathbf Q_n^+\underline w}\rangle.\]
\end{proof}

\subsubsection*{Proof of Theorem \ref{stabilization}}
Let $\bfP\in\Par$ and $n\in\mathbb N$.
Let $S(\bfP_n^+)=\{w\in\widetilde\rW| U_{w(\alpha)}\neq\bfP_n^+
\text{\ for every\ }\alpha\in\widetilde\Delta\}\cup\{1\}$
and $\mathrm Y(\bfP_n^+)=\bigcup_{w\in S(\bfP_n^+)}(LG/\bfI)^{\leq w}\in\Upsilon$.
We claim that for any  
$\mathrm Y\supset\mathrm Y(\bfP_n^+)\in\Upsilon$
we have equalities
\beq\label{equ a}
 \langle A^\mathrm Y_\mM\rangle*\langle\delta_{\bfP_n^+}\rangle=
\langle A^{\mathrm Y(\bfP_n^+)}_\mM\rangle*\langle\delta_{\bfP_n^+}\rangle
\eeq
\beq\label{equ b}
\langle A^\mathrm Y_\mM\rangle*\langle\delta_{\bfP^+}\rangle=\langle\pi_\bfP^!\mM_\bfP\rangle.
\eeq
Since for any $\mF\in M(LG)$ we have 
$\mF\is\delta_{\bfI_n^+}*\mF$ for some $n$, 
the system 
 $\{\langle A^\mathrm Y_\mM\rangle*\langle\mF\rangle\}_{\mathrm Y\in\Upsilon}=
\{\langle A^\mathrm Y_\mM\rangle*\langle\delta_{\bfI_n^+}\rangle*\langle\mF\rangle\}_{\mathrm Y\in\Upsilon}=\{\langle A^{\mathrm Y(\bfI_n^+)}_\mM\rangle*\langle\mF\rangle\}_{\mathrm Y\in\Upsilon}$ stabilizes as $\mathrm Y\supset Y(\bfI_n^+)$.
The theorem follows.

Proof of the claim.
Both equalities~\eqref{equ a} 
and~\eqref{equ b}
were proved
in \cite[Proposition 4.2.4 ]{BKV}
in the case $\mM=\delta$ is the unit.
For the equality~\eqref{equ a},
the proof  uses 
\cite[Lemma 4.2.3]{BKV}
which is our
Lemma \ref{key lemma} in the special case of $\mM=\delta$.
The same argument 
works for general objects in $\frakA(LG)$ once we replace 
Lemma 4.2.3 in \emph{loc. cit.}
by  Lemma \ref{key lemma}.

For the second equality~\eqref{equ b}, the argument in 
\emph{loc. cit.}  implies that (again  replacing 
Lemma 4.2.3 in \emph{loc. cit.}
by  Lemma \ref{key lemma}), for any $\alpha\in\widetilde\Delta$
such that $U_\alpha\in\bfP^+$, we have 
$\langle A^\mathrm Y_\mM\rangle*\langle\delta_{\bfP^+}\rangle=\langle
\pi_{\bfP_{\widetilde\Delta\setminus\alpha}}^!(\mM_{\bfP_{\widetilde\Delta\setminus\alpha}})\rangle*\langle\delta_{\bfP^+}\rangle$. Since $\bfP\subset\bfP_{\widetilde\Delta\setminus\alpha}$,
 we have 
 $\langle\pi_{\bfP_{\widetilde\Delta\setminus\alpha}}^!(\mM_{\bfP_{\widetilde\Delta\setminus\alpha}})\rangle*\langle\delta_{\bfP^+}\rangle\stackrel{~\eqref{convolution 1}}=\langle\pi_\bfP^!(\on{HC}_{\bfP,\bfP_{\widetilde\Delta\setminus\alpha}}(\mM_{\bfP_{\widetilde\Delta\setminus\alpha}}))\rangle=\langle\pi_\bfP^!\mM_\bfP\rangle$
 and hence we conclude that 
 $\langle A^\mathrm Y_\mM\rangle*\langle\delta_{\bfP^+}\rangle=\langle
\pi_\bfP^!\mM_\bfP\rangle$.

\section{Construction of elements in depth zero Bernstein centers}

\subsection{A map from $\frak A(LG)$ to $\frakZ(LG)$}
\label{A to Z}
Giving an object  $\mM\in\frakA(LG)$, 
we define $\langle A_\mM\rangle\in\on{End}(K_0(M(LG)))$
by the formula
\[\langle A_\mM\rangle(\langle\mF\rangle):=\lim_{\mathrm Y\in\Upsilon^{op}}\langle A_\mM^\mathrm Y\rangle*\langle\mF\rangle.\]
Note that the limit exists by the stablization theorem Theorem \ref{stabilization}.
Consider the following algebra 
\[\frakZ(LG):=\on{End}_{K_0(M(LG)^2)}(K_0(M(LG)))\]
which can be regarded as a Grothendieck group version of the Bernstein center.
\begin{thm}\label{geometric Bernstein center}
(1) We have $\langle A_\mM\rangle\in\frakZ(LG)$.

(2) For any $\mathbf{P}\in\on{Par}$,
we have $\langle A_\mM\rangle(\langle\delta_{\mathbf{P}^+}\rangle)=\langle 
\pi_{\mathbf P}^!\mM_{\mathbf P}\rangle$.

(3) 
The assignment 
$\mM\to\langle A_\mM\rangle$ 
satisfies $\langle A_{\mM'}\rangle\circ\langle A_{\mM}\rangle=
\langle A_{\mM'*\mM}\rangle$ 
and hence 
gives rise to an algebra 
homomorphism
\[\langle A\rangle:\frakA(LG)\to\frakZ(LG).\]
\end{thm}
\begin{proof}
Part (2) follows from 
Theorem \ref{stabilization}.
Part (1) was proved in \cite[Theorem 4.1.9]{BKV} in the case
$\mM=\delta$. 
The argument uses  Lemma 4.3.1 and Lemma 4.3.2
in \emph{loc. cit.}, which are general properties about averaging functors on 
$LG$, and hence 
also applies to the general case.

Proof of (3). By the definition of $\langle A_\mM\rangle$, we need to show that 
\[\langle A_{\mM'}\rangle(\langle A^\mathrm Y_{\mM}\rangle)=\langle A^\mathrm Y_{\mM'*\mM}\rangle\]
 for large $\mathrm Y\in\Upsilon$.
 Since $\langle A^\mathrm Y_{\mM}\rangle:=\sum_{\bfP\in\Par}(-1)^{r(G)-r(\bfP)}\langle\on{Av}^{\mathrm Y_\bfP}(h_\bfP^!\mM_\bfP)\rangle$, it will be enough to show, for any $\bfP\in\Par$,
 we have
 \beq\label{mult}
 \langle A_{\mM'}\rangle(\langle\on{Av}^{\mathrm Y_\bfP}(h_\bfP^!\mM_\bfP)\rangle)=\
 \langle\on{Av}^{\mathrm Y_\bfP}(h_\bfP^!(\mM'_\bfP*\mM_\bfP))\rangle.
 \eeq
By additivity, we can assume $\mathrm Y_\bfP=\bfI w\bfP/\bfP$ where 
$w$  is an element in the coset $w\rW_\bfP$ of minimal length. 
Let $\bfI_w=\prod_{\beta\in\tilde{\Phi}^+\cap w^{-1}(\tilde{\Phi}^-)} U_{\beta}$.
Then the projection $LG\to LG/\bfP$ induces an isomorphism 
$\bfI_w\is \mathrm Y_\bfP$ and hence by \cite[Lemma 2.3.2]{BKV}
we have an isomorphism
$\on{Av}^{\mathrm Y_\bfP}(h_\bfP^!\mM_\bfP)\is\on{Av}^{\bfI_w}(\pi_\bfP^!\mM_\bfP)$
for any $\mM_\bfP\in M(\frac{LG/\bfP^+}{\bfP})$. 
On the other hand, the same proof of  \cite[Corollary 4.1.10]{BKV} shows that 
$\langle A_\mM\rangle(\langle\on{Av}^{\bfI_w}(\mF)\rangle)=
\langle\on{Av}^{\bfI_w}\rangle(\langle A_\mM\rangle(\langle\mF\rangle))$
for any $\mF\in M(LG)$. All together, we see that the left hand side of~\eqref{mult}
is equal to
\[\langle A_{\mM'}\rangle(\langle\on{Av}^{\mathrm Y_\bfP}(h_\bfP^!\mM_\bfP)\rangle)=
\langle A_{\mM'}\rangle(\langle\on{Av}^{\bfI_w}(\pi_\bfP^!\mM_\bfP)\rangle)=
\langle\on{Av}^{\bfI_w}\rangle(\langle A_{\mM'}\rangle(\langle\pi_\bfP^!\mM_\bfP\rangle))=\]
\[=\langle\on{Av}^{\bfI_w}\rangle(\langle A_{\mM'}\rangle(\langle\delta_{\bfP^+}\rangle*\langle\pi_\bfP^!\mM_\bfP\rangle))\stackrel{\text{part\ }(2)}=
\langle\on{Av}^{\bfI_w}\rangle(\langle\pi_\bfP^!\mM'_\bfP\rangle*\langle\pi_\bfP^!\mM_\bfP\rangle)=
 \langle\on{Av}^{\mathrm Y_\bfP}(h_\bfP^!(\mM'_\bfP*\mM_\bfP))\rangle.
\]
The desired claim follows. This finishes the proof of (3) and hence the theorem.
  
\end{proof}

\subsection{Depth zero Bernstein centers}\label{depth zero}
In this section consider the case when $k=\overline{\bbF}_q$
and $G$ is defined over $\bbF_q$. 
Let $F=\bbF_q((t))$.
We write $G(F)=LG(\mathbb F_q)$ for the corresponding 
reductive group over the local function field $F$.
Let $\on{Rep}(G(F))$ be the category of smooth representation 
over $\overline\bbQ_\ell$.
Let $Z(G(F))$ be the Bernstein center of $G(F)$, that is,
the $\overline\bbQ_\ell$-algebra of endomorphisms of the 
identity functor $\id_{\on{Rep}(G(F))}$.

Let $(M(G(F)),*)$ be the Hecke algebra of $G(F)$ consisting smooth measures on
$G(F)$ with compact supports. 
Each $z\in Z(G(F))$ defines an endomorphism 
$z_{M(G(F))}$ of the Hecke algebra $M(G(F))$ characterized by the formula:
For every $(\pi,V)\in\on{Rep}(G(F))$, $v\in V$, and $h\in M(G(F))$, we have 
an equality $z(h(v))=(z_{M(G(F))}(h))(v)$.
The action of $G(F)^2$ on $G(F)$ given by 
$(x,y)(g)=xgy^{-1}$ 
defines a $M(G(F))^2$-action on $M(G(F))$ and 
the map 
sending $z$ to $z_{M(G(F))}$ defines an algebra isomorphism 
$Z(G(F))\is\End_{M(G(F))^2}(M(G(F)))$.

The category $\on{Rep}(G(F))$ decomposes into direct sum
$\on{Rep}(G(F))=\on{Rep}^0(G(F))\oplus
\on{Rep}^{>0}(G(F))$ where $\on{Rep}^0(G(F))$ (resp. $\on{Rep}^{>0}(G(F))$)
is the category smooth representations of depth zero (resp. of positive depths).
Thus the Bernstein center $Z(G(F))$ decomposes 
as a direct sum $Z(G(F))=Z^0(G(F))\oplus Z^{>0}(G(F))$
where $Z^0(G(F))$ is the subalgebra of depth zero Bernstein center.

\subsection{Geometric construction of elements in the depth zero Bernstein center}
\label{A to Z}
In this section we preserve the setup in Section \ref{depth zero}.
We assume $G$ is split
and $T$ is a split maximal torus.
Let 
 $\on{Fr}:G\to G$ be the geometric Frobenius. 
 All the geometric objects and categories introduced before are defined over $\mathbb F_q$ and 
we denote by
$M^\Fr(LG)^{}$, $M^\Fr(\frac{LG/\bfP^+}{\bfP})$, etc, the corresponding categories
of 
$\Fr$-equivariant objects.

We  have 
a version of 
Theorem \ref{geometric Bernstein center} 
in the $\Fr$-equivariant setting. Namely,
if we write
\[\frakA^\Fr(LG):=\lim_{\mathbf P\in\Par} K_0(M^\Fr(\frac{LG/\mathbf P^+}{\mathbf P}))\]
\[\frakZ^\Fr(LG):=\on{End}_{K_0((M^\Fr(LG))^2)}(K_0(M^\Fr(LG))).\]
Then for any
 $\mM\in\frakA^\Fr(LG)$ one can associate to it 
an object 
$\langle A^\Fr_{\mM}\rangle\in\frakZ^\Fr(LG)$
such that the assignment 
$\mM\to\langle A_{\mM}^\Fr\rangle$ defines
an algebra homomorphism
\[
\langle A^\Fr\rangle:\frakA^\Fr(LG)\to\frakZ^\Fr(LG)
\ \ \ \ \ \langle A^\Fr\rangle(\langle\mM\rangle)=\langle A_{\mM}^\Fr\rangle\]

Recall the  Bernstein center $Z(G(F))$ of $G(F)$ and 
the subalgebra of depth zero Bernstein center $Z^0(G(F))\subset Z(G(F))$
in Section \ref{depth zero}.
According to the sheaf-function correspondence 
for measures  
developed in \cite[Section 3.4.6]{BKV}, 
we have natural algebra map 
\[K_0(M^{\Fr}(LG))\to M(G(F))\ \ \ \langle \mF\rangle\to [\mF]\]
which induces 
a natural algebra map
\beq\label{End}
\frakZ^\Fr(LG)
\to \End_{M(G(F))^2}(M(G(F)))\is Z(G(F))\ \ \ \ \langle A_{\mM}^\Fr\rangle\to [A_{\mM}^\Fr]
%K_0(A(LG)^\Fr)\to A(G(F)).
\eeq
Thus combining the two above map we obtain an algebra map
\beq\label{formula [A]}
[A^\Fr]:\frakA^\Fr(LG)\to\frakZ^\Fr(LG)\to Z(G(F))\ \ \ [A^\Fr](\mM)=[A_\mM^\Fr]
\eeq

\begin{thm}\label{main}
(1) The image of $[A^\Fr]$ lies in $Z^0(G(F))$ and hence gives rise to an algebra map
\[[A^\Fr]:\frakA(LG)^\Fr\to Z^0(G(F)).\]
(2) For any $\mathbf P\in\Par$, we have 
$[A_\mM^\Fr]([\delta_{\mathbf P^+}])=[\pi^!_\mathbf P\mM_\mathbf P]$.
\end{thm}
\begin{proof}
Let $\mM\in\frakA(LG)^\Fr$.
By (3) Theorem \ref{geometric Bernstein center},
we have 
$[A_\mM^\Fr]=[A_\mM^\Fr* A^\Fr_\delta]=[A_\mM^\Fr]\circ[A_\delta^\Fr]$
where $\delta\in\frakA(LG)^\Fr$ is the unit element.
By \cite[Theorem 4.4.1]{BKV}, the element  $[A_\delta^\Fr]\in Z^0(G(F))$
is the projector to the depth zero spectrum and it implies 
$[A_\mM^\Fr]=[A_\mM^\Fr]\circ[A_\delta^\Fr]\in Z^0(G(F))$. 
Part (2) follows from (2) of
Theorem \ref{geometric Bernstein center}.

\end{proof}

\subsection{The algebra $A(G(F))$}
We shall prove a version of
Theorem \ref{geometric Bernstein center} 
at the level of measures.
For any $\bfP\in\Par$, let  $\mathrm P=\bfP(\bbF_q)$
and $\mathrm P^+=\bfP(\bbF_q)$ be the corresponding 
parahoric subgroup  and  pro-unipotent radical respectively.
We denote by
$(M(\frac{G(F)/\mathrm P^+}{\mathrm P}),*)$ be the parahoric Heck algebra 
of $G(F)$ consisting of $\mathrm P^+$ bi-invaraint
and $\mathrm P$-conjugation invariant   smooth measures on $G(F)$ with compact support.
Consider the following algebra 
\[A(G(F))=\lim_{\bfP\in\Par} M(\frac{G(F)/\mathrm P^+}{\mathrm P})\]
where the limit is taken with respect to the map
$M(\frac{G(F)/\mathrm Q^+}{\mathrm Q})\to M(\frac{G(F)/\mathrm P^+}{\mathrm P})$ sending 
$h$ to $h*\delta_{\mathrm P^+}$, where $\bfP\subset\bfQ\in\Par$ and
$\delta_{\mathrm P^+}$ is the Haar measure of $\mathrm P^+$ with total measure one.

Let $\Omega$ be the set of all finite $\mathrm I$-invariant subset 
$Y\subset G(F)/\mathrm I$ such that 
for all $w'\leq w$ such that $\mathrm Iw\subset Y$ we have 
$\mathrm Iw'\subset Y$. For every $\bfP\in\Par$ and $Y\in\Omega$, we denote
by $Y_\mathrm P\subset G(F)/\mathrm P$ the image of $Y$.
For any measure 
$h_\mathrm P\in M(\frac{G(F)/\mathrm P^+}{\mathrm P})$
and $Y\in\Omega$, we define $\on{Av}^{Y_\mathrm P}(h_\mathrm P)$ to be the measure
\[\on{Av}^{Y_\mathrm P}(h_\mathrm P)=\sum_{y\in Y_\mathrm P}\on{Ad}_{y}(h_\mathrm P)\in M(G(F)).\]
For each $Y\in\Omega$ and $h=\{h_\mathrm P\}_{\bfP\in\Par}\in A(G(F))$,
we set
\beq\label{A^y_h}
[A^Y_h]=\sum_{\bfP\in\Par}(-1)^{r(G)-r(\bfP)}\on{Av}^{Y_\mathrm P}(h_{\mathrm P})\in M(G(F)).
\eeq
We have the following theorem generalizing \cite[version 3, Theorem 5.2.2]{BKV}:

\begin{thm}\label{elementary construction}
(1) For every $f\in M(G(F))$ and $h\in A(G(F))$, the sequence 
$\{[A^\mathrm Y_h]*f\}_{\mathrm Y\in\Omega}$ 
of measures 
stabilizes.
\

(2) For each $h\in A(G(F))$, we define 
$[A_h]\in\End_{M(G(F))^{op}}(M(G(F)))$ by the formula
\[[A_h](f)=\lim_{\mathrm Y\in\Omega}\ [A^Y_h]*f\]
We have $[A_h]\in Z^0(G(F))\subset Z(G(F))=\End_{M(G(F))^2}(M(G(F)))$
and the assignment $h\to [A_h]$ defines an algebra map
\beq\label{[A]}
[A]:A(G(F))\to Z^0(G(F))\ \ \ \ h\to [A](h)=[A_h]
\eeq
\

(3) The map $[A]$ in~\eqref{[A]} fits into the following commutative diagram
\[\xymatrix{
\frakA(LG)^\Fr\ar[d]\ar[r]^{\langle A^\Fr\rangle}&\frakZ(LG)^\Fr\ar[d]\\
A(G(F))\ar[r]^{[A]}&Z^0(G(F))}\]
where the vertical arrows are given by the sheaf-function correspondence 
for measures.
\end{thm}
\begin{proof}
The same arguments of Theorem \ref{stabilization}, 
but omitting all the geometry from here, 
shows that the sequence $\{[A_h^\mathrm Y]*f\}_{\mathrm Y\in\Omega}$
stabilizes. Part (1) follows.
For part (2), we note that  for a fixed $\bfP\in\Par$, if we choose  $\mathrm Y$ to be $\mathrm P$-invariant,
 the measure $[A_h^\mathrm Y]\in M(G(F))$ will be $\mathrm P$-conjugation invariant.
Since $G(F)$ is generated as a group by the standard parahoric subgroups $\mathrm P\subset G(F)$, $\bfP\in\Par$, we conclude that the limit 
$[A_h]\in\End_{M(G(F))^{op}}(M(G(F)))$ is  $G(F)$-conjugation invariant 
and hence $[A_h]\in \End_{M(G(F))^2}(M(G(F)))=Z(G(F))$.
On the other hand, the same arguments of Theorem \ref{A to Z} (again omitting all the geometry) show that 
$[A_{h'}]\circ[A_h]=[A_{h'*h}]$ for all $h',h\in A(G(F))$
and we conclude that the assignment $h\to [A_h]$ defines 
an algebra homomorphism 
$[A]:A(G(F))\to Z(G(F))$. To show that the image of $[A]$ lies in $Z^0(G(F))$
we observe that the image $[A](\delta)=[A_\delta]\in Z^0(G(F))$
of the unit element $\delta=\{\delta_{\mathrm P^+}\}_{\bfP\in\Par}$ 
is the depth zero
spectrum projector  $z^0\in Z^0(G(F))$ 
and hence $[A](h)=[A](h*\delta)=[A_h]\circ [A_\delta]=
[A_h]\circ z^0\in Z^0(G(F))$.
Part (2) follows.
Part (3) follows from the sheaf-function correspondence.

\end{proof}

\begin{remark}\label{mixed char}
The formula~\eqref{A^y_h} makes sense in any non-Archimedean local fields 
and the same arguments  shows that 
 Theorem \ref{elementary construction} remains valid in 
the mixed characteristic case. 
\end{remark}

\section{
Strongly central complexes}\label{SCC}
For the reset of the paper we will preserve the setup in Section \ref{A to Z}.
We will assume chosen an isomorphism 
\beq\label{identification}
\iota:k^\times\is(\mathbb Q/\mathbb Z)_{p'}
\eeq
where $(\mathbb Q/\mathbb Z)_{p'}\subset(\mathbb Q/\mathbb Z)$
is the subgroup consisting of elements of order prime to $p$.

\subsection{Strongly central complexes on torus}
Let $\on{sign}:\rW\to\{\pm1\}$ be the sign character of $\rW$. 
For any tame local system $\mL$ on $T$ we denote by
$\rW'_\mL=\{w\in\rW|w^*\mL\is\mL\}$.
Recall the  definition of strongly central complexes on $T$ in \cite{C1}.
\begin{definition}\label{strongly central complexes}
A $\rW$-equivariant complex $\mF\in D(T/\rW)$ on $T$ is called 
strongly central if for any tame local system $\mL$ on $T$ the natural action of 
$\rW'_\mL$ on $H^*_c(T,\mF\otimes\mL)$ is given by the sign
character $\on{sign}|_{\rW'_\mL}:\rW'_{\mL}\to\{\pm1\}$.
We write $A(T/\rW)\subset D(T/\rW)$ for the subcategory of strongly central complexes on $T$.
\end{definition}

We have the following  key technical result.

\begin{thm}\label{BK}
Let $\mF\in A(T/\rW)$ be a strongly central complex.
Then
$\on{HC}_{P,G}(\Ind_{T\subset B}^G(\mF)^\rW)$ is supported on $\frac{P/U_P}{P}$
for all standard parabolic subgroups $P\in\on{par}$.

\end{thm}
\begin{proof}
The case when $P=B$ is the Borel subgroup 
is the main result in \cite{C2}. 
The general case follows from Corollary \ref{B implies P}.
\end{proof}

\subsection{Examples}\label{example}
We recall a construction of strongly central local systems 
in \cite[Section 5.1]{C2}.
Let $\pi_1^t(T)$ be the tame fundamental group of $T$
and let 
$\pi_1(T)_\ell$ be its  pro-$\ell$ quotient.
Let $R_T=\on{Sym}(\pi_1(T)_\ell\otimes_\bbZ\overline\bbQ_\ell)$
the symmetric algebra of $\pi_1(T)_\ell\otimes_{\bbZ_\ell}\overline\bbQ_\ell$
and $R_T^+$ be its argumentation ideal.
For any continuous character
$\chi:\pi_1^t(T)\to\overline\bbQ_\ell^\times$ 
we denote by 
$\rW'_\chi$ the stabilizer of $\chi$ in $\rW$.
The group $\rW_\chi$ acts on $R_T$ and 
$R_T^+$ and the quotient 
$R_\chi=R_T/\langle R_{T,+}^{\rW'_\chi}\rangle$
(where $\langle R_{T,+}^{\rW'_\chi}\rangle$ is the ideal generated by 
the $\rW'_\chi$-invariants
$R_{T,+}^{\rW'_\chi}$) is naturally a 
$\pi_1^t(T)\rtimes\rW'_\chi$-representation and hence gives rise to a 
$\rW_\chi$-equivariant $\overline\bbQ_\ell$-adic tame local system 
$\mE_\chi^{uni}$ on $T$.
Since $\mL_\chi$ is $\rW'_\chi$-equivariant, the 
 tensor product  $\mE_\chi=\mE^{uni}_\chi\otimes\mL_\chi$
is also $\rW'_\chi$-equivariant.
The
induced
$\rW$-equivariant local system 
$\Ind_{\rW'_\chi}^\rW\mE_\chi$ 
depends only on the
$\rW$-orbit
$\theta=\rW\chi$ of $\chi$
 and we denote by
$\mE_\theta=\Ind_{\rW'_\chi}^\rW\mE_\chi$ the resulting $\rW$-equivariant local system
and $\mE_\theta^\vee$ the dual of $\mE_\theta$.

\begin{lemma}\label{E_theta}
The 
$\rW$-equivariant local system 
$\mE_\theta^\vee$ is a strongly central.
\end{lemma}
\begin{proof}
Let $\mL$ be a tame local system on $T$.
Then,
up to degree shifts,
there is a $\rW_\mL$-equivariant isomorphism between
$H_c^*(T,\mE_\theta^\vee\otimes\mL^{-1})$ and $H^*(T,\mE_\theta\otimes\mL)$
 and hence it suffices to 
show that $\rW_\mL$-acts on $H^*(T,\mE_\theta\otimes\mL)$ 
via the sign character. 
Recall the $\ell$-adic Mellin transform
$\mathfrak M: D(T)\is D_{Coh}^b(\calC(T))$
where $\calC(T)$ is the $\overline\bbQ_\ell$-scheme of tame local systems on $T$
(see, e.g., \cite[Section 4.1]{C1}).
There is a $\rW'_\mL$-equivariant isomorphism 
$H^*(T,\mE_\theta\otimes\mL)\is Li_\mL^*\frakM(\mE_\theta)$
where $i_\mL:\{\mL\}\to\calC(T)$ is the inclusion.
According to 
\cite[Lemma 4.3]{C2}, the restriction of 
$\frakM(\mE_\theta\otimes\on{sign})$ 
to the component 
$\calC(T)_\mL\subset\calC(T)$ containing $\mL$ descends to
the quotient 
$\calC(T)_\mL//\rW'_\mL$ and it follows that 
the action of $\rW'_\mL$ on the derived fiber 
$Li_\mL^*\frakM(\mE_\theta\otimes\on{sign})$ is trivial. The desired claim follows.\footnote{In \cite{C2},  we work with 
a certain subgroup 
 $\rW_\mL\subset\rW_\mL'$  (see \cite[Section 4.2]{C2})
 instead of the full stabilizer subgroup $\rW'_\mL$ (we have $\rW_\mL=\rW_\mL'$ if the center of $G$ is connected).
However, the same argument in \emph{loc. cit.} goes through if 
we used $\rW_\mL'$.}

\end{proof}

In  \cite{BK2}, Braverman-Kazhdan associated to each representation 
$\rho:\hat G\to\GL_n$ of the dual group and a non-trivial character 
$\psi:\bbF_q\to\overline\bbQ_\ell^\times$, 
a $\rW$-equivariant perverse sheaf $\mF_{T,\rho,\psi}$ on $T$, called the 
$\rho$-Bessel sheaf, with remarkable properties (see Section \ref{Deligne} for the relationship between $\mF_{T,\rho,\psi}$
and Deligne's $\epsilon$-factors).\footnote{The Bessel sheaf
$\mF_{T,\rho,\psi}$  was denoted by $\Phi_{T,\rho}$ in \emph{loc. cit.}.}
In \cite[Theorem 1.4]{C2} we showed 
that $\rho$-Bessel sheaf 
$\mF_{T,\rho,\psi}$ is strongly central.
Now
Theorem \ref{BK} implies the following conjecture of Braverman-Kazhdan:

\begin{corollary}\cite[Conjecture 6.5.]{BK2}\label{BK conj}
The Harish-Chandra transfrom
$\on{HC}_{P,G}(\Ind_{T\subset B}^G(\mF_{T,\rho,\psi})^\rW)$ is supported on $\frac{P/U_P}{P}$
for all standard parabolic subgroups $P\in\on{par}$.
\end{corollary}

\subsection{From strongly central complexes to $\frakA(LG)_1$}\label{admissible to strongly central complex}
For any $\bfP\subset\bfQ\in\Par$,
the pull-back along 
the natural identification 
$\phi_{\bfP,\bfQ}:T_\bfP\is T_\bfQ$ induces 
a functor 
$\phi^!_{\bfP,\bfQ}:D(T_\bfQ/\rW_\bfQ)\to D(T_\bfP/\rW_\bfP)$.
It follows from the definition  of  strongly central complexes
that $\phi^!_{\bfP,\bfQ}$ restricts to a functor 
\[\phi^!_{\bfP,\bfQ}:A(T_\bfQ/\rW_\bfQ)\to A(T_\bfP/\rW_\bfP)\]
and hence we can form the following limit
\[\on{lim}_{\bfP\in\Par}K_0(A(T_\bfP/\rW_\bfP))\]
where the limit is taken with respect to $\phi^!_{\bfP,\bfQ}$.
 For any $\bfP\in\Par$
and $\mF_{\bfP}\in D(T_\bfP/\rW_\bfP)$  we write
\beq\label{F_L}
\mF_{L_\bfP}:=\Ind_{T_\bfP\subset B_\bfP}^{L_\bfP}(\mF_{T_\bfP})^{\rW_\bfP}
[-2\dim U_\bfP](-\dim U_\bfP)\in D(\frac{L_\bfP}{L_\bfP}).
\eeq
where $\mF_{T_\bfP}\in D(\frac{T_\bfP}{T_\bfP})$ is the $!$-pullback of 
$\mF_\bfP$ long the projection $\frac{T_\bfP}{T_\bfP}\to T_\bfP/\rW_\bfP$.

\begin{lemma}\label{T/W to G}
(1) The assignment $\{\langle\mF_{T_\mathbf P}\rangle\}_{\mathbf P\in\Par}\to\{\langle\mF_{L_\mathbf P}\rangle\}_{\mathbf P\in\Par}$.
defines a natural  map
\beq\label{lim to A}
\on{lim}_{\mathbf P\in\Par}K_0(A(T_\mathbf P/\rW_\mathbf P))\to \on{lim}_{\mathbf P\in\Par}K_0(A(L_\mathbf P))=\frakA(LG)_1
\eeq
(2)
The assignment $\langle\mF\rangle\to \{\langle\phi_\mathbf P^!\mF\rangle\}_{\mathbf P\in\Par} $ (where $\phi_\mathbf P:T_\mathbf P\is T$ is the natural identification)
defines a map
\beq\label{section T/W}
\eta:K_0(A(T/\rW))\to \lim_{\mathbf P\in\Par} K_0(A(T_\mathbf P/\rW_\mathbf P))
\eeq

\end{lemma}
\begin{proof}
By Lemma \ref{W summand} and Theorem \ref{BK}, we have $\mF_{L_\bfP}\in A(L_\bfP)$
and, for any $\bfP\subset\bfQ\in\Par$, we have 
\[\on{Res}_{L_\bfP\subset B_{\bfP,\bfQ}}^{L_\bfQ}{\mF_{L_\bfQ}}[2\dim U_{\bfP,\bfQ}](\dim U_{\bfP,\bfQ})\is
\on{Res}_{L_\bfP\subset B_{\bfP,\bfQ}}^{L_\bfQ}
\Ind_{T_\bfQ\subset B_\bfQ}^{L_\bfQ}(\mF_{T_\bfQ})^{\rW_\bfQ}[-2\dim U_\bfP](-\dim U_\bfP)\is\]
\[\is\Ind_{T_\bfP\subset B_\bfP}^{L_\bfP}(\mF_{T_\bfP})^{\rW_\bfP}
[-2\dim U_\bfP](-\dim U_\bfP)\is
\mF_{L_\bfP}\]
(note that $\dim U_{\bfP,\bfQ}=\dim U_\bfQ-\dim U_\bfP$).
The proposition follows.

\end{proof}

Combining with Proposition \ref{1=par} and Lemma \ref{T/W to G}, we obtain 
a map
\beq\label{limit to admissible}
 \Phi:\lim_{\bfP\in\Par} K_0(A(T_\bfP/\rW_\bfP))\to \lim_{\bfP\in\Par} K_0(A(L_\bfP))=\frakA(LG)_1\subset\frakA(LG)
\eeq
In Section \ref{de-cate}  we shall relate the map above 
with a certain map from the limit of stable Bernstein center of
finite reductive groups to the depth zero Bernstein center.

\subsection{Deligne-Lusztig packets}\label{DL-packets}

Let $\on{Irr}(G(\bbF_q))$ be the set of isomorphism classes of irreducible representations 
of the finite reductive group $G(\bbF_q)$ over $\overline{\bbQ}_\ell$.
Let $\hat G$ be the Langlands dual group of $G$ over $\overline{\bbQ}_\ell$.
In \cite{DL}, Deligne and Lusztig 
proved the following results
\begin{thm}\label{DL}
(1) There is natural a natural bijection
\footnote{ 
The bijection depends on the isomorphism~\eqref{identification}}
between the set of 
$G(\bbF_q)$-conjugacy classes of  pairs $(S,\chi)$, 
where $S$ is a $\Fr$-stable maximal torus of $G$ and 
$\chi:S(\bbF_q)\to\overline{\bbQ}_\ell^\times$ is a character,
and the set of semi-simple conjugacy classes in 
$\hat G$ stable under $x\to x^q$.
For any such pair $(S,\chi)$ we denote by $\theta$ the corresponding
semi-simple conjugacy classe in $\hat G$.
(2) For any irreducible representation $\pi\in\on{Irr}(G(\bbF_q))$, there exists 
a pair $(S,\chi)$ as above such that 
$\langle\pi,R_{S,\chi}\rangle\neq 0$, where $R_{S,\chi}$ is the Deligne-Lusztig 
virtual representation associated to $(S,\chi)$. Moreover, the 
semi-simple conjugacy class $\theta$ 
of $(S,\chi)$  is uniquely determined by 
$\pi$.
\end{thm}

Let $\hat T\subset\hat G$ be the canonical split maximal torus.
The map $q:\hat T\to \hat T, t\to t^q$
is $\rW$-equivariant and hence descends to
a map on the adjoint quotient 
$[q]:\hat T//\rW\to\hat T//\rW$. Moreover, there is a natural bijection between
the set of semi-simple conjugacy classes in 
$\hat G$ stable under $x\to x^q$ and the  
 the fixed point set $(\hat T//\rW)^{[q]}$.
Thus the theorem above gives rise to 
a well-define surjective map
\beq\label{DL-map}
\frak{L}:\on{Irr}(G(\bbF_q))\to(\hat T//\rW)^{[q]}
\eeq
sending $\pi$ to the corresponding semi-simple conjugacy class $\theta\in(\hat T//\rW)^{[q]}$.
The fibers $\frak{L}^{-1}(\theta)$ of the map~\eqref{DL-map}
are called the  Deligne-Lusztig packets.

\subsection{Stable Bernstein center for finite reductive groups
and functoriality}\label{center for finite gps}
We study stable Bernstein center
for $G(\bbF_q)$.

We first recall the definition of Bernstein center
for $G(\bbF_q)$.
Let $\on{Rep}(G(\bbF_q))$ be the category of finite dimensional representation 
of $G(\bbF_q)$ over $\overline\bbQ_\ell$.
Let $Z(G(\bbF_q))=\on{End}_{}(\Id_{\on{Rep}(G(\bbF_q))})$
be the Bernstein center of $\on{Rep}(G(\bbF_q))$, that is, the endomorphism 
of the identity functor $\Id_{\on{Rep}(G(\bbF_q))}$.
By definition $Z(G(\bbF_q))$ is a commutative algebra over $\overline\bbQ_\ell$.

Each $z\in Z(G(\bbF_q))$ defines a function 
$\gamma_z:\on{Irr}(G(\bbF_q))\to\overline\bbQ_\ell$, to be called the 
$\gamma$-function of $z$, 
charactered by the formula
$z|_\pi=\gamma_z(\pi)\Id_\pi$ for any $\pi\in\on{Irr}(G(\bbF_q))$.

Each $z\in Z(G(\bbF_q))$ gives rise to a $\overline\bbQ_\ell$-valued
class function
$\beta_z:G(\bbF_q)\to\overline\bbQ_\ell$ charactered by the formula
$z|_\pi=\sum_{g\in G(\bbF_q)} \beta_z(g)\pi(g)$. Moreover 
the assignment 
$z\to \beta_z$ defines an isomorphism of algebras 
$Z(G(\bbF_q))\is C(G(\bbF_q))$, where 
$C(G(\bbF_q))$ is the algebra of $\overline\bbQ_\ell$-valued
classes function with algebra structure given by 
convolution 
$\beta_1*\beta_2(g)=\sum_{h\in G(\bbF_q)}\beta_1(h)\beta_2(h^{-1}g)$.

An element $z\in Z(G(\bbF_q))$ is called \emph{stable} 
if the corresponding $\gamma$-function 
$\gamma_z:\on{Irr}(G(\bbF_q))\to\overline\bbQ_\ell$ 
is constant on Deligne-Lusztig packets, that is,
we have $\gamma_z(\pi)=\gamma_z(\pi')$
if $\pi,\pi'\in\frak L^{-1}(\theta)$ for some $\theta\in(\check T//\rW)^{[q]}$
(see Section \ref{DL-packets}).
The stable Bernstein center for $G(\bbF_q)$ is defined as 
$Z^{st}(G(\bbF_q)):=\{z\in Z^{}(G(\bbF_q))| z\text{\ is stable}\}$.
By definition, $Z^{st}(G(\bbF_q))$ is a subalgebra of $Z(G(\bbF_q))$.

A class function $\beta\in C(G(\bbF_q))$
is called \emph{stable} if it lies in the image of 
$Z^{st}(G(\bbF_q))$ under the isomorphism 
$Z(G(\bbF_q))\is C(G(\bbF_q))$.
 We write $C^{st}(G(\bbF_q))\subset C(G(\bbF_q))$
 for the subalgebra of stable class functions.

 Note that  each $z\in Z^{st}(G(\bbF_q))$ corresponds to
a unique function $f_z:(\check T//\rW)^{[q]}\to\overline\bbQ_\ell$
characterized by the formula
$\gamma_z=f_z\circ\frak L$, where $\frak L$ is the map in~\eqref{DL-map}.
Moreover, the assignment 
$z\to f_z$ defines an algebra isomorphism 
\beq\label{Z=T//W}
Z^{st}(G(\bbF_q))\is\overline\bbQ_\ell[(\hat T//\rW)^{[q]}]
\eeq
where right hand side is the algebra of $\overline\bbQ_\ell$-valued functions 
on the (finite set) $(\check T//\rW)^{[q]}$.

For any $\bfP\in\Par$, 
the canonical identification $\phi_\bfP:T\is T_\bfP$
gives rise to map 
$\hat\phi_\bfP:\hat T_\bfP\is\hat T$ compatible with 
the Weyl group action $\rW_\bfP$ (where $\rW_\bfP$-acts  
on $T$ via the map 
$\rW_\bfP\to\widetilde\rW\to\rW$).
The map $\hat\phi_\bfP$ induces a map
\beq\label{[q]}
\hat\phi_\bfP^{[q]}:(\hat T_\bfP//\rW_\bfP)^{[q]}\to(\hat T//\rW)^{[q]}
\eeq
and we define 
$
\hat\rho_\bfP:Z^{st}(G(\bbF_q))\to Z^{st}(L_\bfP(\bbF_q))
$
to be the following composition
\beq\label{G to L_P}
\hat\rho_\bfP:Z^{st}(G(\bbF_q))\is\overline\bbQ_\ell[(\hat T//\rW)^{[q]}]\stackrel{(\hat\phi_\bfP^{[q]})^*}\longrightarrow
\overline\bbQ_\ell[(\hat T_\bfP//\rW_\bfP)^{[q]}]\is Z^{st}(L_\bfP(\bbF_q))
\eeq
where the first and last maps are the isomorphisms in~\eqref{Z=T//W}.

Similarly, for any $\bfP\subset\bfQ\in\Par$,
the natural $\rW_\bfQ$-equivariant 
identification 
$\hat\phi_{\bfP,\bfQ}:\hat T_\bfP\is \hat T_\bfQ$ induces a 
natural map
\beq\label{[q]_P,Q}
\hat\phi_{\bfP,\bfQ}^{[q]}:(\hat T_\bfP//\rW_\bfP)^{[q]}\to(\hat T_\bfQ//\rW_\bfQ)^{[q]}
\eeq
which gives rise to a map
\beq\label{transfer Q to P}
\hat\rho_{\bfP,\bfQ}:Z^{st}(L_\bfQ(\bbF_q))\to Z^{st}(L_\bfP(\bbF_q))
\eeq
By construction, for any $\bfP\subset\bfQ\in\Par$
we have
\beq\label{comp 1}
\hat\rho_{\bfP,\bfQ}\circ\hat\rho_{\bfQ}=\hat\rho_{\bfP}
\eeq
and for any 
$\bfP\subset\bfQ\subset\bfQ'\in\Par$
we have 
\beq\label{comp 2}
\hat\rho_{\bfP,\bfQ}\circ\hat\rho_{\bfQ,\bfQ'}=\hat\rho_{\bfP,\bfQ'}
\eeq
Introduce the following spaces 
 \beq
 \on{colim}_{\bfP\in\Par}(\hat T_\bfP//\rW_\bfP)^{[q]},\ \ \ \ \
\lim_{\bfP\in\Par}Z^{st}(L_\bfP(\bbF_q)) \eeq
where the colimit is taken with respect to the maps
 $\hat\phi_{\bfP,\bfQ}^{[q]}$ and  the limit is taken with respect to $\hat\rho_{\bfP,\bfQ}$.
 We have 
\beq\label{lim=colim}
\lim_{\bfP\in\Par} Z^{st}(L_\bfP(\bbF_q))\is 
 \lim_{\bfP\in\Par}\overline\bbQ_\ell[(\hat T_\bfP//\rW_\bfP)^{[q]}]\is 
\overline\bbQ_\ell[\on{colim}_{\bfP\in\Par}(\hat T_\bfP//\rW_\bfP)^{[q]}].
\eeq

The map
~\eqref{comp 1} induces a
natural map 
\beq\label{section}
 i:\on{colim}_{\bfP\in\Par}(\hat T_\bfP//\rW_\bfP)^{[q]}\to (\hat T//\rW)^{[q]}
\eeq
and we have a commutative diagram
\beq\label{hat rho}
\xymatrix{Z^{st}(G(\bbF_q))\ar[r]^{\on{lim}_{}\hat\rho_\bfP\ \ \ \ \ }\ar[d]_{~\eqref{Z=T//W}}^\simeq&\lim_{\bfP\in\Par}Z^{st}(L_\bfP(\bbF_q))\ar[d]_{~\eqref{lim=colim}}^\simeq\\
\overline\bbQ_\ell[T//\rW)^{[q]}]\ar[r]^{i^*\ \ \ \ \ \ \ }&\overline\bbQ_\ell[\on{colim}_{\bfP\in\Par}(\hat T_\bfP//\rW_\bfP)^{[q]}]}.
\eeq

\subsection{From stable center for finite reductive groups to 
the depth zero Bernstein center}\label{de-cate}
For any   
$\Fr$-equivariant strongly central complex 
$\mF\in A^\Fr(T/\rW)$, the associated
complex
$\underline\mF_G=\underline\Ind_{T\subset B}^G(\underline\mF)^\rW\in D^\Fr(G)$
is naturally $\Fr$-equivariant  
and we let
\[\beta_\mF=
\on{Tr}(\Fr,\underline\Ind_{T\subset B}^G(\underline\mF)^\rW)\in C(G(\bbF_q))\]
be the corresponding class function on $G(\bbF_q)$ and
$z_\mF\in Z(G(\bbF_q))$ be the corresponding element 
in the Bernstein center.
  \begin{prop}\label{geometrization 1}
  For any $\mF\in A^\Fr(T/\rW)$  we have 
$z_\mF\in Z^{st}(G(\bbF_q))$
and the assignment $\mF\to z_\mF$
induces a surjective map
\beq\label{surjection}
\upsilon:K_0(A^\Fr(T/\rW))\to Z^{st}(G(\bbF_q)).
\eeq
\end{prop}
 \begin{proof}
 To prove that $z_\mF$ is stable,
 it is equivalent to show that the corresponding class function $\beta_\mF$ is stable.
 For this it suffices to show that 
 for any pair $(S,\chi)$ as in Theorem \ref{DL} we have 
 \beq\label{stability}
 \beta_\mF*\Tr(R_{S,\chi})=\gamma\cdot\Tr(R_{S,\chi})
 \eeq
 where $\Tr(R_{S,\chi})$ is the character of the Deligne-Lusztig virtual representation $R_{S,\chi}$ and
 $\gamma\in\overline{\mathbb Q}_\ell$ is a constant.
We will use the following result of Lusztig.
Recall that 
 for split group $G$ over $k$ the
$G(\bbF_q)$-conjugacy classes of $\Fr$-stable maximal torus 
 are in bijection with conjugacy classes of the Weyl group $\rW$.
Let $w\in\rW$ be a representative of the 
Weyl group conjugacy class corresponding to $S$.
Then each character  
 $\chi$ of $S(\bbF_q)$ corresponds to 
a character local system $\mL_\chi$ on $T$
such that $\Fr^*\mL_\chi\is (w^{-1})^*\mL_\chi$. 
A choice of an isomorphism above 
endows the induction $\underline\Ind_{T\subset B}^G(\mL_\chi)$ 
a canonical $\Fr$-equivariant structure and   Lusztig 
\cite[Proposition 8.15 and 9.2]{Lu1} (see also \cite[Theorem 3.7]{BK2})
proved that 
\[\Tr(R_{S,\chi})=q^{\dim G/B}\Tr(\Fr,\underline\Ind_{T\subset B}^G(\mL_\chi)).\]
Thus to prove~\eqref{stability}, it suffices to 
construct a $\Fr$-equivariant isomorphism
\beq\label{action}
\underline\mF_G
*\underline\Ind_{T\subset B}^G(\mL_\chi)\is \on{H}_c^*(T,\underline\mF\otimes\mL^{-1}_\chi)\otimes
\underline\Ind_{T\subset B}^G(\mL_\chi)
\eeq
(note that the $\rW$-equivariant structure of $\mF$ together with the
isomorphism $\Fr^*\mL_\chi\is (w^{-1})^*\mL_\chi$ endows the
the vector space $\on{H}_c^*(T,\underline\mF\otimes\mL^{-1}_\chi)$  a 
$\Fr$-action). 
By Lemma \ref{basic compatibility} and Corollary \ref{B implies P} we have a canonical  isomorphism
\[\underline\mF_G*\underline\Ind_{T\subset B}^G(\mL_\chi)\is\underline\Ind_{T\subset B}^G(\underline\Res_{T\subset B}^G(\underline\mF_G)*\mL_\chi)\]
On the other hand,  Lemma \ref{W summand} implies
\[\underline\Res_{T\subset B}^G(\underline\mF_G)*\mL_\theta\is\underline\mF*\mL_\chi
\is\on{H}_c^*(T,\underline\mF\otimes\mL^{-1}_\chi)\otimes\mL_\chi\]
Combining the above two isomorphisms we obtain the  isomorphism in~\eqref{action}
and 
one can check that it is compatible with the natural 
$\Fr$-equivariant structure. The desired claim follows.

We shall show that the map $\upsilon$ is surjective.
Let $\theta$ be a $\rW$-orbit of a tame character 
and let 
$\mE_\theta^\vee\in A(T/\rW)$ be the strongly central local system in Lemma \ref{E_theta}.
Assume $\theta$ is Frobenius-invariant. Then 
$\mE_\theta^\vee\in
A^\Fr(T/\rW)$ 
and  by~\eqref{action} we have 
\[\underline\mF_{G,\theta}*\underline\Ind_{T\subset B}^G(\mL_\chi)\is\on{H}_c^*(T,\underline\mE_\theta^\vee\otimes\mL_\chi^{-1})\otimes\underline\Ind_{T\subset B}^G(\mL_\chi)\]
where $\underline\mF_{G,\theta}=\underline\Ind_{T\subset B}^G(\underline\mE_\theta^\vee)^\rW$.
Note that  
$\on{H}_c^*(T,\underline\mE_\theta^\vee\otimes\mL_\chi^{-1})=0$
if $\chi^{-1}\notin\theta$
and $\on{H}_c^*(T,\underline\mE_\theta^\vee\otimes\mL_\chi^{-1})\is\on{H}^*(T,\underline\mE_\chi^{uni})$
if $\chi^{-1}\in\theta$ (here $\mL_\chi$ is the local system above).
Thus by taking trace of Frobenius of the above isomorphism, we obtain 
\beq\label{quasi-idempotent}
\beta_{\mE_{\theta}^\vee}*\Tr(R_{S,\chi})=\gamma\cdot\Tr(R_{S,\chi})
\eeq
where $\beta_{\mE_{\theta}^\vee}=\Tr(\Fr,\underline\mF_{G,\theta})$
and the constant $\gamma$ is given by:
$\gamma=0$ if $\chi^{-1}\notin\theta$
and $\gamma=\Tr(\Fr,\on{H}^*(T,\mE_\chi^{uni}))$ if 
$\theta=\rW\chi^{-1}$.
Note that the number
\beq\label{gamma_theta}
\gamma_\theta:=\Tr(\Fr,\on{H}^*(T,\underline\mE_{\chi^{-1}}^{uni}))=\Tr(\Fr,\on{H}^*(T,\underline\mE_\chi^{uni}))
\eeq
depends only on the $\rW$-orbit $\theta=\rW\chi^{-1}$ (note that $\mE_\chi^{uni}=\mE_{\chi^{-1}}^{uni}$) 
and we claim that $\gamma_{\theta}\neq 0$.
Then
the image 
$\upsilon(\langle\mE_{\theta}^\vee\rangle)$
is the non-zero function given by 
\beq\label{char function}
\upsilon(\langle\mE_{\theta}^\vee\rangle)=\gamma_\theta\cdot 1_{\theta^{-1}}
\eeq
where $1_{\theta^{-1}}$ is the characteristic function supported at $\theta^{-1}:=\rW\chi\in(\hat T//\rW)^{[q]}$.
Since there is a bijection between 
$(\hat T//\rW)^{[q]}$
and the set of Frobenius-invariant
$\rW$-conjugacy classes of tame characters of $T$, 
the collection $\{\gamma_\theta\cdot 1_{\theta^{-1}}\}_{\theta\in(\hat T//\rW)^{[q]}}$ form a basis of 
$\overline\bbQ_\ell[(\hat T//\rW)^{[q]}]$. It 
implies $\upsilon$ is surjective.

Proof of the claim.
Since the $\overline\bbQ_\ell$-Tate module 
$\pi_1(T)_\ell\otimes_{\bbZ_\ell}\overline\bbQ_\ell\is\overline\bbQ_\ell(1)^{\oplus\dim T}$
is of weight $-2$
and $\mE_{\chi}^{un}$ is the unipotent local system corresponds to the 
$R_T=\on{Sym}(\pi_1(T)_\ell\otimes_{\bbZ_\ell}\overline\bbQ_\ell)$-module 
$R_T/\langle R_{T,+}^{\rW_\chi}\rangle$, it follows that 
there is filtration of $\mE_{\chi}^{un}$ with associated grade
$\on{gr}(\mE_{\chi}^{un})=\bigoplus_{i=0}^m\overline\bbQ_\ell(i)^{\oplus n_i}$,
$m$ and $n_i\in\bbZ_{\geq 0}$
and it implies $\gamma_\theta=\bigoplus_{i= 0}^m n_i\Tr(\Fr,\on{H}^*(T,\overline\bbQ_\ell(i)^{}))=
(\oplus_{i=0}^m n_iq^{-i})\cdot\Tr(\Fr,\on{H}^*(T,\overline\bbQ_\ell))\neq 0$.

 \end{proof}

\begin{prop}\label{geometrization 2}
For any $\mF=\{\langle\mF_{\mathbf P}\rangle\}_{\mathbf P\in\Par}\in\lim_{\mathbf P\in\Par} K_0(A^\Fr(T_\mathbf P/\rW_\bfP))$ the collection
$z_\mF=\{z_{\mF_{\mathbf P}}\}_{\mathbf P\in\Par}$ defines an element in
$\lim_{\mathbf P\in\Par} Z^{st}(L_\mathbf P(\bbF_q))$ and the assignment 
$\mF\to z_\mF$ induces a surjective 
map
\beq\label{surjection lim}
\Upsilon:\lim_{\mathbf P\in\Par} K_0(A(T_\mathbf P/\rW_\mathbf P)^\Fr)\to \lim_{\mathbf P\in\Par} Z^{st}(L_\mathbf P(\bbF_q)).
\eeq
%Moreover, we have the following diagram commutes
%\[\xymatrix{K_0(A(T/\rW)^\Fr)\ar[r]^{\eta\ \ \ \ }\ar[d]^\upsilon&\lim_{\mathbf P\in\Par} K_0(
%A(T_\mathbf P/\rW_\mathbf P)^\Fr)\ar[d]^\Upsilon\\
%Z^{st}(G(\bbF_q))\ar[r]^{i\ \ }&
 %\lim_{\mathbf P\in\Par}Z^{st}(L_\mathbf P(\bbF_q))}\]

\end{prop}
 \begin{proof}
Let  $\bfP\subset\bfQ\in\Par$ and
let $\pi$ and $\pi'$ be a representation of  $L_\bfP(\bbF_q)$ and $L_\bfQ(\bbF_q)$
with Deligne-Lusztig 
 semi-simple parameter $\theta_\bfP\in(\hat T_\bfP//\rW_\bfP)^{[q]}$ 
 and $\theta_{\bfQ}=\hat\rho_{\bfP,\bfQ}^{[q]}(\theta_\bfP)\in(\hat T_\bfQ//\rW_\bfQ)^{[q]}$ respectively. 
It follows from the isomorphism~\eqref{action} that 
\[z_{\mF_\bfP}|_\pi=\Tr(\Fr,\on{H}_c^*(T_\bfP,\underline\mF_\bfP\otimes\mL_\chi^{-1})\cdot\id_\pi\ \ \ \ z_{\mF_\bfQ}|_{\pi'}=\Tr(\Fr,\on{H}_c^*(T_\bfQ,\underline\mF_\bfQ\otimes\mL_\chi^{-1})\cdot\id_{\pi'}\]
where $\chi$ (resp. $\chi'$)
is any tame character of $T_\bfP$ (resp. $T_\bfQ$)
for which the corresponding 
semi-simple conjugacy class in $\hat L_\bfP$ (resp. $\hat L_\bfQ$) is equal to 
$\theta_\bfP$ (resp. $\theta_\bfQ$).
To show that the map $\Upsilon$ is well defined, it suffices to show that 
$\Tr(\Fr,\on{H}_c^*(T_\bfP,\underline\mF_\bfP\otimes\mL_\chi^{-1})=
\Tr(\Fr,\on{H}_c^*(T_\bfQ,\underline\mF_\bfQ\otimes\mL_{\chi'}^{-1})$
where $\chi$ and $\chi'$ above satisfy
 $\chi=\chi'\circ\rho_{\bfP,\bfQ}$.  
This follows from the fact that the pull-back along the $\rW_\bfP$-equivariant isomorphism $\rho_{\bfP,\bfQ}:T_\bfP\is T_\bfQ$ induces a $\Fr$-equivariant isomorphism
\[\on{H}_c^*(T_\bfQ,\underline\mF_\bfQ\otimes\mL_{\chi'}^{-1})\is
\on{H}_c^*(T_\bfP,\rho_{\bfP,\bfQ}^!\underline\mF_\bfQ\otimes\rho_{\bfP,\bfQ}^!\mL_\chi^{-1})\is
\on{H}_c^*(T_\bfP,\underline\mF_\bfP\otimes\mL_{\chi}^{-1}).\]
 
We shall show that $\Upsilon$ is surjective. 
Let $\theta=\rW\chi^{-1}\in(\hat T//\rW)^{[q]}$ and let
$s_{\theta}=i^{-1}(\theta^{-1})\subset\on{colim}_{\bfP\in\Par}(\hat T_\bfP//\rW_\bfP)^{[q]}$
the pre-image of $\theta^{-1}=\rW\chi$ along the map in
~\eqref{section}. 
Recall the local system 
$\mE_{\theta}^\vee\in A^\Fr(T/\rW)$ in Lemma \ref{E_theta} attached to 
$\theta$.
According to Lemma \ref{T/W to G}, the collection  
$\eta(\mE_{\theta}^\vee)=\{\phi_\bfP^!(\mE_{\theta}^\vee)\}_{\bfP\in\Par}$
defines an element in $\lim_{\bfP\in\Par} A^\Fr(T_\bfP/\rW_\bfP)$.
Moreover, it follows from~\eqref{char function} that
the element
$\Upsilon(\langle\eta(\mE_\theta^\vee)\rangle)\in\lim_{\bfP\in\Par} Z^{st}(L_\bfP(\bbF_q))\is\overline\bbQ_\ell[\on{colim}_{\bfP\in\Par}(\hat T_\bfP//\rW_\bfP)^{[q]}]$
is equal to 
\beq\label{char s_theta}
\Upsilon(\langle\eta(\mE_\theta^\vee)\rangle)=\gamma_\theta\cdot1_{s_{\theta}}
\eeq
where $\gamma_\theta$ is the non-zero constant in~\eqref{gamma_theta} and 
$1_{s_\theta}$ the characteristic function of the subset $s_\theta$.

Let $e_\bfP:(\hat T_\bfP//\rW_\bfP)^{[q]}\to\on{colim}_{\bfP\in\Par}(\hat T_\bfP//\rW_\bfP)^{[q]}$
be the natural map and, for any point 
$s_{\theta,i}\in s_\theta=\{s_{\theta,1},...,s_{\theta,l_\theta}\}$,
we denote by
$\theta_{\bfP,i}=e_\bfP^{-1}(s_{\theta,i})$.
Then the local system 
$\phi_\bfP^!(\mE_\theta^\vee)$ admits the following direct sum decomposition
\beq\label{mon decomp}
\phi_\bfP^!(\mE_\theta^\vee)=\oplus_{i=1}^j\phi_\bfP^!(\mE_\theta^\vee)_i
\eeq
where $\phi_\bfP^!(\mE_\theta^\vee)_i$ is the summand whose $\ell$-adic Mellin transform
is set theoretically supported on $\theta_{\bfP,i}$ (here we view $\theta_{\bfP,i}$
as a collection of tame characters on $T_\bfP$ and hence 
a subset of $\overline\bbQ_\ell$-scheme of
 tame local systems 
$\calC(T_\bfP)$ on $T_\bfP$).
Note that we have $\phi_\bfP^!(\mE_\theta^\vee)_i\in A^\Fr(T_\bfP/\rW_\bfP)$
and 
the collection 
$\eta(\mE_\theta^\vee)_i:=\{\phi_\bfP^!(\mE_\theta^\vee)_i\}_{\bfP\in\Par}$
defines an element in $\on{colim} A^\Fr(T_\bfP/\rW_\bfP)$.
Moreover, the decomposition~\eqref{mon decomp} induces an isomorphism 
\[\eta(\mE_\theta^\vee)=\oplus_{i=1}^j \eta(\mE_\theta^\vee)_i\]
and~\eqref{char s_theta} implies 
\[\Upsilon(\langle\eta(\mE_\theta^\vee)_i\rangle)=\gamma_\theta\cdot 1_{s_{\theta,i}}.\]
where $1_{s_{\theta,i}}$ is the characteristic function of the point $s_{\theta,i}$.
Since the collection $\{1_{s_{\theta,i}}\}$, $\theta\in(\hat T//\rW)^{[q]}$ 
and $i\in[1,l_\theta]$ forms a basis of $\overline\bbQ_\ell[\on{colim}_{\bfP\in\Par}(\hat T_\bfP//\rW)^{[q]}]$
we conclude that $\Upsilon$ is surjective.

\end{proof}

 \begin{thm}\label{decategorification}
We have the following commutative diagram 
\beq\label{key diagram}
\xymatrix{K_0(A^\Fr(T/\rW))\ar[r]^{\eta^\Fr\ \ \ \ \ }\ar[d]^{\upsilon}&
\lim_{\mathbf P\in\Par} 
K_0(A^\Fr(T_\mathbf P/\rW_\mathbf P))\ar[r]^{\ \ \ \ \ \ \ \ \ \Phi^\Fr}\ar[d]^{\Upsilon}&\frakA^\Fr(LG)\ar[d]^{[A^\Fr]}\\
Z^{st}(G(\bbF_q))\ar[r]^{\lim\hat\rho_\bfP\ \ \ \ }&
 \lim_{\mathbf P\in\Par}Z^{st}(L_\mathbf P(\bbF_q))\ar[r]^{\ \ \ \ \ \Psi}&Z^0(G(F))}
 \eeq
where $\Phi^\Fr$ is the map in~\eqref{limit to admissible}
and the map $\Psi$ is given by the formula:
for any element $z=\{z_\mathbf P\}\in\lim_{\mathbf P\in\Par}Z^{st}(L_\mathbf P(\bbF_q))$,
$(\pi,V)\in\on{Rep}^0(G(F))$ and a vector 
$v\in V$
we have 
\beq\label{formula}
(\Psi(z)|_\pi)(v)=(z_\mathbf P|_{\pi^{P^+}})(v)
\eeq
where $\mathbf P\in\Par$ is any parahoric subgroup such that 
$v\in V^{P^+}$
and $\pi^{P^+}$ denotes the natural representation of 
$L_\mathbf P(\bbF_q)$ on $V^{P^+}$.
 \end{thm}
 \begin{proof}
Consider the following diagram commutes
 \[\xymatrix{\bfP\ar[r]^{\pi_\bfP}\ar[d]_{q_\bfP}&\frac{\bfP}{\bfP}\ar[d]^{f_\bfP}\\
 L_\bfP\ar[r]^{\pi_{L_\bfP}}&\frac{L_\bfP}{L_\bfP}}\]
 Let $\mF=\{\mF_\bfP\}\in \lim_{\mathbf P\in\Par} K_0(A^\Fr(T_\bfP/\rW_\mathbf P))$
 and let $z=\Upsilon(\mF)$ and $\mM=\Phi^\Fr(\mF)$.
By  Proposition \ref{1=par} and Lemma \ref{geometrization 2} we have 
\[\mM=\{\langle\mM_\bfP\rangle=\langle f_\bfP^!\mF_{L_\bfP}\rangle\}_{\bfP\in\Par}\in\frakA^\Fr(LG)_1\subset\frakA^\Fr(LG) \]
and
\[z=\{z_{\bfP}=z_{\mF_\bfP}\}_{\bfP\in\Par}\in\lim_{\mathbf P\in\Par}Z^{st}(L_\mathbf P(\bbF_q))\]
According to Theorem \ref{main}, for any 
$\bfP\in\Par$, we have 
 \[([A^\Fr](\mM))([\delta_\bfP^+])=[A_\mM^\Fr]([\delta_\bfP^+])=[\pi_\bfP^!f_\bfP^!\mF_{L_\bfP}]
=[q_\bfP^!\pi_{L_\bfP}^!\mF_{L_\bfP}]\]
  Note that, by Lemma \ref{basic compatibility} we have
 \[\pi_{L_\bfP}^!\mF_{L_\bfP}\is\pi_{L_\bfP}^!\Ind_{T_\bfP\subset B_\bfP}^{L_\bfP}(\mF_{T_\bfP})^{\rW_\bfP}
 [-2\dim U_\bfP](-\dim U_\bfP)\is
\underline\Ind_{T_\bfP\subset B_\bfP}^{L_\bfP}(\underline\mF_\bfP)^{\rW_\bfP}
\] 
Now Lemma \ref{geometrization 1} implies that 
 for any depth zero representation $(\pi,V)$ and each $v\in V^{P^+}$  
 we have 
 \[([A^\Fr](\mM))(v)=([A^\Fr](\mM))(\delta_{\mathrm P^+}(v))=
 (([A^\Fr](\mM))([\delta_\bfP^+]))(v)=
 [q_\bfP^!\pi_{L_\bfP}^!\mF_{L_\bfP}](v)=\]
 \[=\sum_{g\in L_\bfP(\bbF_q)}\beta_{\mF_{\bfP}}(x)\pi^{P^+}(x)(v)=(z_{\mF_\bfP}|_{\pi^{P^+}})(v)=(\Psi(z)|_\pi)(v).
 \]
 The theorem follows.
 \end{proof}

 \section{Applications}
  We discuss  
 geometric construction of Delinge-Lusztig parameters for depth zero representations
 and certain remarkable elements in 
$Z^0(G(F))$ coming from Deligne's epsilon factors \cite{D}.

\subsection{Deligne-Lusztig parameters}\label{DL parameters}
Theorem \ref{decategorification}  gives rise to an embedding
\beq\label{zeta}
\zeta:=\Psi\circ\lim \hat\rho_\bfP:Z^{st}(G(\bbF_q))\to\lim_{\bfP\in\Par}
Z^{st}(L_\bfP(\bbF_q))
\to Z^0(G(F))
\eeq
Since there is a natural bijection between 
the set of characters $\on{Hom}(Z^{st}(G(\bbF_q)),\overline\bbQ_\ell)$ 
with
 the set $(\hat T//\rW)^{[q]}$ of 
 semi-simple conjugacy classes in the dual group $\hat G$ (over $\overline\bbQ_\ell$) stable under $x\to x^q$,
one can associate to each
irreducible depth zero representation $\pi$ of $G(F)$ 
a point $\theta(\pi)$ in $(\hat T//\rW)^{[q]}$, to be called the 
Deligne-Lusztig parameter of $\pi$,
 obtained by
composing
$\zeta$
 with the 
central character $Z^0(G(F))\to\on{End}(\pi)=\overline\bbQ_\ell$.
Using formula~\eqref{formula}, one can check that the 
Deligne-Lusztig parameter of $\theta(\pi)$ of $\pi$ agrees with the one 
coming from Moy-Prasad theory \cite{MP1,MP2}.

 \subsection{Bernstein centers arising from Deligne's epsilon factors}\label{Deligne}
Let $W_F$ be the Weil group of $F$. We fix a geometric Frobenius 
element $\Fr\in W_F$
and denote by $v:W_F\to\bZ$ the canonical map
sending $\Fr$ to $1$.
Let $I\subset\rW_F$ be the inertia group and $P\subset I$ be the 
wild inertia subgroup. 

A Langlands parameter 
over $\overline\bbQ_\ell$ is a pair
$(r,N)$ where 
$r:W_F\to \hat G(\overline\bbQ_\ell)$
is a continuous homomorphism with open kernel such that 
$r(w)$ is semisimple for all $w\in W_F$ and  
$N\in\hat\fg(\overline\bbQ_\ell)$ is a nilpotent element 
such that 
$\on{Ad}_{\rho(w)}(N)=q^{-v(w)}N$. A Langlands parameter $r$ is called tame if 
$P\subset\on{ker}r$.
We call
two parameters  $(r,N)$ and $(r',N')$ equivalent if 
there is a $x\in \hat G(\overline\bbQ_\ell)$
such that $\on{Ad}_x(r(w))=r'(w)$ for all 
$x\in W_F$ and $\on{Ad}_x(N)=N'$.
We denote by $\on{Loc}_{\hat G,F}$ and  $\on{Loc}^t_{\hat G,F}$
the sets of equivalence classes of Langlands parameters 
and  tame Langlands parameters respectively.

Inspired by the work of  Macdonald \cite{M}, 
we call two  Langlands parameters 
$(r,N)$ and $(r',N')$  $I$-equivalent 
if the restrictions $r_I$ and $r_I'$ of
$r$ and $r'$ to the inertia group $I$ are equivalent, that is,
there exists an element $x\in \hat G(\overline\bbQ_\ell)$
such that $\on{Ad}_x(r(w))=r'(w)$ for all 
$x\in I$. Note that the $I$-equivalence classes of 
$(r,N)$ only depends on $r$ and 
if $(r,N)$ and $(r',N')$ are $I$-equivalent 
and $(r,N)$ is tame then $(r',N')$ is also 
tame.
We denote by 
$\Phi_{\hat G,F,I}$ and 
$\Phi_{\hat G,F,I}^t$ the set of 
$I$-equivalence classes of Langlands parameters
and 
tame Langlands parameters.

\quash{
%Since 
%$(r,N)$ is $I$-equivalent to $(r,0)$
%(here $0\in\hat\fg(\overline\bbQ_\ell)$ is the zero element), 
there is a bijection between 
$\on{Loc}_{\hat G,F,I}$ (resp. 
$\on{Loc}_{\hat G,F,I}^t$)
and the the $I$-equivalence classes of 
Langlands parameters (resp. tame Langlands parameters)
of the form $(r,0)$.
}

Let  $H$ be a split reductive group with dual group $\hat H$ over 
$\overline\bbQ_\ell$. Then any 
group homomorphism 
$\rho:\hat 
H\to \hat G$ between dual groups induces 
a map 
$\rho:\Phi_{\hat H,F}\to\Phi_{\hat G,F}$,
$\rho^t:\Phi_{\hat H,F}^t\to\Phi_{\hat G,F}^t$,
and 
$\rho^t_I:\Phi_{\hat H,F,I}^t\to\Phi_{\hat G,F,I}^t$.

\begin{lemma}\label{transfer}
The isomorphism $\iota:k^\times\is(\mathbb Q/\mathbb Z)_{p'}$
in~\eqref{identification}
gives rise to a bijection
$\Phi_{\hat G,F,I}^t\stackrel{\simeq}\to (\hat T//\rW)^{[q]}$.\footnote{The author learned this fact from C.C. Tsai.}
Moreover, for any $\rho:\hat H\to \hat G$ as above, we have the following comutative diagram
\[\xymatrix{\Phi_{\hat H,F,I}^t\ar[r]^{\rho^t_I}\ar[d]^{\simeq}&\Phi_{\hat G,F,I}^t\ar[d]^{\simeq}\\
(\hat T_H//\rW_H)^{[q]}\ar[r]^\rho&(\hat T//\rW)^{[q]}}.\]
\end{lemma}
\begin{proof}
Let $(r,N)$ be a tame Langlands parameter.
Since $P\subset\on{ker}(r)$ the restriction $r_I=r|_I$
of 
$r$ to $I$ factors through the quotient
$\bar r_I:I/P\to\hat G(\overline\bbQ_\ell)$.
Since $I/P\is\lim_{n\in\bbZ_{>0}}\mu_n(k)\is\Hom_{\bbZ}((\mathbb Q/\mathbb Z)_{p'},k^\times)$,  the isomorphism $\iota:(\mathbb Q/\mathbb Z)_{p'}\is k^\times$ gives rise to 
a pro generator $\sigma\in I/P$ and we let
$s=\bar r_I(\sigma)\in\hat G(\overline\bbQ_\ell)$.
Since $\Fr\sigma^q\Fr^{-1}=\sigma$,  the two elements
$s$ and $s^q$ are conjugated in $\hat G(\overline\bbQ_\ell)$ 
and the map 
sending $(r,N)$ to $s$ defines a bijection 
between the set $\Phi_{\hat G,F,I}^t$ of $I$-equivalent classes of tame
Langlands
parameters and the set $(\hat T//\rW)^{[q]}$ of semisimple conjugacy classes in 
$\hat G(\overline\bbQ_\ell)$ stable under 
the map $x\to x^q$. The second claim is clear.
\end{proof}

\begin{remark}
The above bijection was first
observed in \cite{M} for $G=\GL_n$.
\end{remark}

Consider the case when $G=\GL_n$ is the general linear group.
Let $\psi_F:F\to\overline\bbQ_\ell^\times$ be a nontrivial character of 
$F$ such that $\fp\subset\on{ker}\psi_F$ 
but $\mO_F\not\subset\on{ker}\psi_F$.
The restriction of $\psi_F$ to $\mO_F$ descends to a non-trivial character 
$\psi:\bbF_q\is\mO_F/\fp_F\to\overline\bbQ_\ell^\times$.
Let 
$dx$ be the Haar measure on $F$ such that 
$\fp_F$ has mass $q^{-1/2}$.
In \cite[Theorem 6.5]{D}, Deligne associated to 
each $(r,N)\in\on{Loc}_{\GL_n,F}$ a nonzero constant
\beq
\epsilon_0(r,\psi_F,dx)\in\overline\bbQ_\ell^\times
\eeq
to be called the epsilon factor for $(r,N)$.
When $(r,N)\in\Phi_{\GL_n,F}^t$, the constant 
is calculated in \cite[Section 5]{D} and it follows that 
$\epsilon_0(r,\psi,dx)$ depends only on the 
$I$-equivalence class of $(r,N)$ and hence gives rise to a function 
\beq\label{deligne}
\epsilon_{0}:\Phi_{\GL_n,F,I}^t\to\overline\bbQ_\ell^\times\ \ \ \  \epsilon_{0}((r,N))=\epsilon_0(r,\psi_F,dx).
\eeq

Let now $G$ be any split reductive group.
Assume that we are given  
a representation $\rho:\hat G\to\GL_n$ of the dual group $\hat G$.
Then the pullback of  $\epsilon_{0}$
along the map 
$\rho^t_I:\Phi_{\hat G,F,I}^t\to\Phi^t_{\GL_n,F,I}$
defines an element:
\[
\epsilon_{0,\rho}:=\epsilon_0\circ\rho^t_I:\Phi_{\hat G,F,I}^t\to\Phi_{\GL_n,F,I}^t\to\overline\bbQ_\ell^\times.
\]
Using Lemma \ref{transfer}, one can view 
$\epsilon_{0,\rho}$ as 
an element in $Z^{st}(G(\bbF_q))\is\overline\bbQ_\ell[(\hat T//\rW)^{[q]}]\is\overline\bbQ_\ell[\Phi_{\hat G,F,I}^t]$ 
and we denote by
\beq\label{z_rho}
z_{0,\rho}:=\zeta(\epsilon_{0,\rho})\in Z^0(G(F))
\eeq
the image of $\epsilon_{0,\rho}$ under the embedding 
$\zeta:Z^{st}(G(\bbF_q))\to Z^0(G(F))$
in~\eqref{zeta}.
Now
a conjecture of Braverman and Kazhdan in \cite{BK2}, proved in \cite{C1}, 
imply the following geometric formula for $z_{0,\rho}$. 
Assume that $G$ is  semi-simple and simply connected.
Recall the Braverman-Kazhdan 
Bessel sheaf  
$\mF_{T,\check\rho,\psi}\in A(T/\rW)$ on $T$ associated to 
the dual representation $\check\rho$ and $\psi$  
in Section \ref{example}. 
The Bessel sheaf $\mF_{T,\check\rho,\psi}$ has a canonical $\Fr$-equivariant 
structure coming from the one on the Artin-Schreier sheaf $\mL_\psi$
and we denote by 
\[z_{\check\rho}=[A^\Fr]\circ\Phi^\Fr\circ\eta^\Fr(\langle\mF_{T,\check\rho,\psi}\rangle)\in Z^0(G(F))\]
the image under the composed map
$[A^\Fr]\circ\Phi^\Fr\circ\eta^\Fr:K(A^\Fr(T/\rW))\to\lim_{\bfP\in\Par} K_0(A^\Fr(T_\bfP/\rW_\bfP))\to\frakA^\Fr(LG)\to Z^0(G(F))$ in Theorem \ref{decategorification}.
\begin{thm}\label{deligne}
We have 
$z_{\check\rho}=(-1)^nz_{0,\rho}$.
\end{thm}
\begin{proof}
We will view
$\epsilon_{0}$ and $\epsilon_{0,\rho}$
as elements in $Z^{st}(\GL_n(\bbF_q))$ and 
$Z^{st}(G(\bbF_q))$ via the identifications
$\overline\bbQ_\ell[\Phi_{\GL_n,F,I}^t]\is \overline\bbQ_\ell[(\hat T_n//\rW_n)^{[q]}]\is Z^{st}(\GL_n(\bbF_q))$ and
$\overline\bbQ_\ell[\Phi_{\hat G,F,I}^t]\is \overline\bbQ_\ell[(\hat T//\rW)^{[q]}]\is Z^{st}(G(\bbF_q))$.
A result of Macdonald \cite{M} says that 
$\epsilon_{0}|_\pi=(-1)^n\gamma_{\psi}(\check\pi)\id$, where 
$\gamma_{\psi}:\on{Irr}(\GL_n(\bbF_q))\to\overline\bbQ_\ell$ 
is the $\gamma$-function for $\GL_n(\bbF_q)$
(see line (2.5) of \cite{M} or line (1.3) of \cite{BK2}
for the definition of $\gamma_\psi$).\footnote{In \cite[line (2.5)]{M}, 
the $\gamma$-function is denoted by
$\gamma_\psi(\pi)=w(\pi,\psi)$.}
In view of Lemma \ref{transfer}, we see that 
$\epsilon_{0,\rho}|_\pi=(-1)^n\gamma_{G,\rho,\psi}(\check\pi)\id$
where $\gamma_{G,\rho,\psi}:\on{Irr}(G(\bbF_q))\to\overline\bbQ_\ell$
is the Braverman-Kazhdan $\gamma$-function in \cite[Section 1.4]{BK2}.
On the other hand,
the proof of Braverman-Kazhdan conjecture in \cite[Corollary 1.7]{C1}
implies that 
$\upsilon(\langle\mF_{T,\check\rho,\psi}\rangle)|_\pi=\gamma_{G,\rho,\psi}(\check\pi)\id$ 
(where 
$\upsilon$ is the map in Lemma \ref{surjection})
and hence we conclude that \[\epsilon_{0,\rho}=(-1)^n\upsilon(\langle\mF_{T,\check\rho,\psi}\rangle)\in Z^{st}(G(\bbF_q)).\]
Since the element $z_{0,\rho}$
is equal to the image 
$z_{0,\rho}=\zeta(\epsilon_{0,\rho})=\Psi\circ \lim\hat\rho_\bfP(\epsilon_{0,\rho})$ along the bottom arrows in diagram~\eqref{key diagram}, the commutativity of 
the diagram in \emph{loc. cit.}
implies 
\[z_{\check\rho}=[A^\Fr]\circ\Phi^\Fr\circ\eta^\Fr(\langle\mF_{T,\check\rho,\psi}\rangle)=
\Psi\circ \lim\hat\rho_\bfP\circ\upsilon(\langle\mF_{T,\check\rho,\psi}\rangle)
=
(-1)^n\zeta(\epsilon_{0,\rho})=(-1)^nz_{0,\rho}.\]

\end{proof}


\begin{thebibliography}{99999}


\bibitem[BK1]{BK1}A. Braverman, D. Kazhdan, $\gamma$-functions of representations
and lifting (with an appendix by V. Vologodsky), Geom. Funct.
Anal., Special Volume (2000), Part I, 237-278.

\bibitem[BK2]{BK2} A. Braverman, D. Kazhdan, $\gamma$-sheaves on reductive groups,
Studies in Memory of
Issai Schur, Chevaleret-Rehovot, 2000, Birkhäuser Boston, Boston, MA, 2003, 27-47.

\bibitem[BKV]{BKV}
R. Bezrukavnikov, D. Kazhdan, Y. Varshavsky, 
A categorical approach to the stable center conjecture,
Ast\'erisque No. 369 (2015), 27-97


\bibitem[BT]{BT}
R. Bezrukavnikov, K. Tolmachov, 
Exterior powers of a parabolic Springer sheaf, 
arXiv:2209.01603



\bibitem[BV]{BV}
R. Bezrukavnikov, Y. Varshavsky, Affine Springer fibers and depth zero L-packets, arXiv:2104.13123.

\bibitem[C1]{C1} 
T.-H. Chen, A vanishing conjecture: $\on{GL}_n$ case,
Sel. Math. New Ser. 28, 13 (2022). 

\bibitem[C2]{C2} 
T.-H. Chen, On a conjecture of Braverman-Kazhdan,
J. Amer. Math. Soc. 35 (2022), 1171-1214

\bibitem[C3]{C3}
T.-H. Chen,
Towards the depth zero stable Bernstein center conjecture II: Semisimple Langlands parameters,
in preparation. 

\bibitem[C4]{C4}
T.-H. Chen,
Towards the depth zero stable Bernstein center conjecture III: Stability.
in preparation. 

\bibitem[CN]{CN} S. Cheng, B.C. Ng\^o, On a conjecture of Braverman and Kazhdan, 
IMRN, 20
(2018), 6177-6200.

\bibitem[D]{D}
P. Deligne, Les constantes des equations functionnelles des functions L,
Proc. Antwerpen Conference, vol. 2; Lecture Notes in Mathematics
349 (1973), 501-597

\bibitem[DL]{DL} 
P. Deligne, G. Lusztig, Representations of reductive groups over finite fields, Annals of Math. 103 (1976), 103-161.


\bibitem[FS]{FS} 
L. Fargues, P. Scholze, Geometrization of the local Langlands correspondence, 
arXiv:2102.13459.



\bibitem[G]{G} V. Ginzburg, 
Admissible modules on symmetric spaces, 
Ast\'erisque 173-174 (1989), 199-256.

\bibitem[H]{H} T. Haines, 
The stable Bernstein center and test functions for Shimura varieties, 
Automorphic Forms and Galois Representations, London Mathematical Society Lecture Note Series, 118-186.

\bibitem[HZ]{HZ} T. Hemo, X. Zhu, in preparation.

\bibitem[Lu1]{Lu1} G. Lusztig, Character sheaves I, 
Advances in Mathematics
Volume 56, Issue 3, June 1985, 193-237

\bibitem[Lu2]{Lu2} G. Lusztig, Classification of unipotent representations of simple p-adic groups,
Int. Math. Res. Not. (1995), 517-589.


\bibitem[LG]{LG}V. Lafforgue, A. Genestier,
Chtoucas restreints pour les groupes reductifs et parametrisation de Langlands locale, arXiv:1709.00978

\bibitem[LL]{LL} 
G. Laumon E. Letellier, Note on a conjecture of Braverman-Kazhdan,
arXiv:1906.07476.

\bibitem[M]{M}
I. G. Macdonald, Zeta functions attached to finite general linear groups,
Math. Ann. 249, 1-15 (1980). 

\bibitem[MP1]{MP1}
A. Moy, G. Prasad, Unrefined minimal K-types for p-adic groups, Inventiones Mathematicae, 116 (1), 393-408.

\bibitem[MP2]{MP2}
A. Moy, M. Prasad, Jacquet functors and unrefined minimal K-types, Commentarii Mathematici Helvetici. European Mathematical Society Publishing House. 71 (1), 98-121.






\bibitem[SS]{SS} 
P. Scholze, S.W. Shin,
On the cohomology of compact unitary group Shimura varieties at ramified split places, J. Amer. Math. Soc. 26 (2013), no. 1, 261-294.

\bibitem[Y]{Y} Z.Yun,
The spherical part of local and global Springer actions, Math Annalen, 359 (2014), no. 3-4, 557-594.

\bibitem[Z1]{Z1} X. Zhu, 
Affine Grassmannians and the geometric Satake in mixed characteristic,
Annals of Mathematics 185 (2017), 403-492

\bibitem[Z2]{Z2} X. Zhu, 
Coherent sheaves on the stack of Langlands parameters, arXiv:2008.02998.

\end{thebibliography}
\end{document}